\documentclass[12pt]{article}
\scrollmode

\usepackage{amsfonts}
\usepackage{amsmath}
\usepackage{amssymb}
\usepackage{graphicx}

\allowdisplaybreaks

\def\eps{\varepsilon}

\def\P{{\mathbb P}}
\def\E{{\mathbb E}}
\def\Z{{\mathbb Z}}
\def\N{{\mathbb N}}
\def\R{{\mathbb R}}
\def\Q{{\mathbb Q}}

\def\1{{\mathbf 1}}

\def\a{\alpha}
\def\o{\omega}  \def\oo{{\tilde{\omega}}}
\def\l{\lambda}
\def\g{\gamma}

\def\e{\epsilon}
\def\d{\delta}
\def\leps{\ell_\e} \def\O{\Omega}

\def\n{{\cal N}}
\def\w{{\cal W}}
\def\m{{\cal M}}
\def\L{{\cal L}}
\def\k{{\cal K}}

\def\j{{\cal J}}
\def\ce{{\cal E}}
\def\B{{\cal B}}

\def\S{{\cal S}}
\def\C{{\cal C}}

\def\x1{X_1^{(1)}}

\def\nn{\nonumber}

\newcommand{\stack}[2]{\genfrac{}{}{0pt}{3}{#1}{#2}}

\newtheorem{theo}{Theorem}[section]
\newtheorem{prop}[theo]{Proposition}
\newtheorem{lm}[theo]{Lemma}

\newtheorem{rmk}[theo]{Remark}
\newtheorem{df}[theo]{Definition}
\def\beq{\begin{equation}}
\def\eeq{\end{equation}}
\newcommand{\bei}{\begin{itemize}}
\newcommand{\eei}{\end{itemize}}
\newcommand{\ben}{\begin{enumerate}}
\newcommand{\een}{\end{enumerate}}
\newcommand{\beqn}{\begin{eqnarray}}
\newcommand{\beqnn}{\begin{eqnarray*}}
\newcommand{\eeqn}{\end{eqnarray}}
\newcommand{\eeqnn}{\end{eqnarray*}}
\newcommand{\brm}{\begin{rmk}}
\newcommand{\erm}{\end{rmk}}

\title{On symmetric random walks with 
       random conductances on  $\Z^d$
}
\author{L.~R.~G.~Fontes \footnote{IME-USP. Rua do Mat\~ao 1010, 05508-090,
S\~ao Paulo, SP,  Brazil, lrenato@ime.usp.br}
\and 
P.~Mathieu \footnote{CMI, 39 rue Joliot-Curie, 13013 Marseille, FRANCE, 
pierre.mathieu@cmi.univ-mrs.fr}      }

\begin{document}
\maketitle

\begin{abstract} 
We study models of continuous time, symmetric, $\Z^d$-valued random walks in random 
environments. One of our aims is to derive estimates on the decay of 
transition probabilities in a case where a uniform ellipticity assumption is 
absent. 
We consider the case  of  independent conductances with a polynomial tail near 
$0$ and obtain precise asymptotics for the annealed return probability 
and convergence times for the random walk confined to a finite box.  
\end{abstract}



\section{Introduction}

We study continuous time, irreducible, symmetric, nearest neighbor random walks in random 
environments on $\Z^d$. Our aim is to derive estimates on the decay of 
transition probabilities in the absence of a uniform ellipticity assumption.  

The paper has four sections (other than this introduction). Sections~\ref{sec:comp}
and~\ref{sec:dec} deal with the decay of the mean or annealed return probability. 
In Section~\ref{sec:comp}, we consider quite general reversible random walks
in a random environment and  
we establish a comparison lemma for the 
annealed return probability. The proof is based on a trace formula 
(in fact an extension of the trace formula for central probability for 
random walks on amenable groups, see~\cite{kn:PS}).  
In  Section~\ref{sec:dec}, we derive sharp bounds on the decay of the annealed
return probability from direct investigation of traces and eigenvalues when the
rates are i.i.d. random variables chosen from a law with polynomial tail near $0$.  
We then prove that one might
get the classical $t^{-d/2}$ decay or a slower decay of the form $t^{-\g}$, 
where $\g<d/2$ is related to the tail of the law of the rates near $0$. 
In Section~\ref{sec:quenc} we deal with the quenched decay and obtain a partial
result (Theorem~\ref{t'}) that nonetheless establishes a difference with respect
to the annealed decay for small values of $\g$.

In Section~\ref{sec:times}, we discuss finite volume random walks taking their
values in a torus. 
We obtain some quenched estimates on 
convergence times when the random rates are i.i.d., chosen 
from a law with polynomial tail near $0$. 
These  follow from sharp bounds on the spectral gap. In particular we prove a universal lower 
bound for the spectral gap of a symmetric random walk on a torus of side length
$N$ (Proposition~\ref{prop:bound} below)
which allows to separate the effects of the usual diffusive $N^{-2}$ factor
and the contribution of small values of the rates. 

The paper is written in such a way to ease independent reading of the
different parts at the cost of some repetition. 
Sections~\ref{sec:comp} and~\ref{sec:times} are
self-contained; only the spectral gap lower bound~(\ref{eq:bound}) from Section~\ref{sec:times}
is needed to proceed through Section~\ref{sec:dec}.


\vskip 1cm
\section{A comparison lemma for the annealed return probability}
\setcounter{equation}{0}
\label{sec:comp}

We study a family of symmetric, irreducible, nearest neighbors Markov chains 
taking their values in $\Z^d$ 
and constructed in the following way. Let $\O$ be the set of functions 
$\o:\Z^d\times\Z^d\to\R_+$ such that 
$\o(x,y)>0$ iff $x\sim y$, and $\o(x,y)=\o(y,x)$. 
($y\sim x$ means that $x$ and $y$ are nearest neighbors.) 
We call elements of $\O$ environments. 

Define  the Markov generator
\begin{equation}
\label{gen0}
{\cal G}^\o f(x)=\sum_{y\sim x}\o(x,y)\,[f(y)-f(x)].
\end{equation}

As usual,  $\{X_t,\,t\in\R_+\}$ will be the coordinate process on path space 
$(\Z^d)^{\R_+}$ and we use the notation  $\P_x^\o$ to denote the unique probability measure 
on path space under which $\{X_t,\,t\in\R_+\}$ is the Markov process 
generated by~(\ref{gen0}) and satisfying $X_0=x$. Under $\P_x^\o$, $X_0=x$; then the process waits 
for an exponentially distributed random time of parameter $\sum_{y\sim x}\o(x,y)$ and jumps to 
point $x_1$ with probability $\o(x,x_1)/ \sum_{y\sim x}\o(x,y)$; this
procedure is then iterated 
choosing independent hoping times. 
Equivalently, one can define $\P_x^\o$ using the theory of symmetric Dirichlet forms, see~\cite{kn:FOT}.  
The reference space is then $L^2(\Z^d)$, equipped with the counting measure. For functions $f$ and $g$ 
with finite support, let 
\beqnn 
{\cal D}^\o(f,g)=\frac 12 \sum_{x\sim y\in\Z^d} \o(x,y)\,[f(x)-f(y)]\,[g(x)-g(y)]. 
\eeqnn
The bilinear form ${\cal D}^\o$ is closable and its closure is a regular, symmetric Dirichlet form. 
Thus, there exists a Hunt process associated to ${\cal D}^\o$. Note that points have non zero capacity. 
Therefore, the measure $\P_x^\o$ is uniquely determined by ${\cal D}^\o$. It is easy to prove that 
both constructions yield the same law $\P_x^\o$. 

Since $\o(x,y)>0$ for all  neighboring pairs $(x,y)$, $X_t$ is irreducible under $\P_x^\o$ 
for all $x$.
The counting measure on $\Z^d$ is reversible because we have assumed that $\o(x,y)=\o(y,x)$. 

We now choose the rates $\o$ at random, according to a translation invariant law $\Q$ on $\O$. 

In the sequel $\Q.\P_x^\o$ will be used as a short hand notation for the annealed 
law defined by $\Q.\P_x^\o[\,\cdot\,]=\int P_x^\o[\,\cdot\,]\,d\Q(\o)$. Note that $X_t$ is Markov under 
$\P^\o_x$ for any $\o$, but is not Markov anymore under $\Q.\P_x^\o$
for nontrivial $\Q$. 
Let $\P^\o=\P_0^\o$ and $\Q.\P^\o=\Q.\P_0^\o$. 

We are interested in estimating the decay of the annealed return probability 
$\Q.\P^\o[X_t=0]$, as $t$ tends to $+\infty$.  

As a subset of $(\R_+)^{\Z^d\times\Z^d }$, $\O$ is a partially ordered set. By duality, 
one can define a partial order on the set of probabilities on $\O$
in the following way. Given two probabilities, 
$\Q$ and $\Q'$, we say that $\Q'\geq \Q$ if, for any measurable, bounded, increasing  function 
$f:\O\to\R$, we have $\Q'(f)\geq \Q(f)$. ($f$ is increasing if, whenever $\o,\o'\in\O$ 
satisfy $\o'(x,y)\geq\o(x,y)$ for all $x,y$, then $f(\o')\geq f(\o)$.) 

\brm The function $\o\rightarrow \P^\o[X_t=0]$ is not monotonous in $\o$. It
is clearly not increasing. It is also not very difficult to find subgraphs of  
$\Z^d$ for which the removal of an edge decreases the value of
$\P^\o[X_t=0]$  
(left as an exercice), 
which implies that the function $\o\rightarrow
\P^\o[X_t=0]$ 
is not decreasing.   
\erm

\begin{lm}
\label{compalemma} 
\label{comparison}
Let $\Q$ and $\Q'$ be two probabilities on $\O$ such that $\Q'\geq \Q$. 
Assume that for $\Q'+\Q$-almost all $\o$, the Markov chain $X_t$ is conservative under 
$\P^\o$. Then, for all time $t$, we have  
\beq\nn
\Q'.\P^\o[X_t=0]\leq \Q.\P^\o[X_t=0].
\eeq
\end{lm}

\noindent{\bf Proof} We prove that $\Q.\P^\o[X_t=0]$ can be written as a supremum of the 
$\Q$-expectation of decreasing in $\o$ functions. More precisely, let
$B_N=[-N,N]^d$ be the box centered at 
the origin and of radius $N$. Let 
${\cal G}^{\o,N}$ be the restriction of the operator ${\cal G}^\o$ to $B_N$ with Dirichlet boundary conditions 
outside $B_N$ (that is, ${\cal G}^{\o,N}$ is the generator of the process
which coincides with the one given by ${\cal G}^\o$ until the latter process 
leaves $B_N$ for the first time, and then it is killed). 
Then $-{\cal G}^{\o,N}$ is a positive symmetric operator. 
Let $\{\mu_i^\o(B_N), i\in[1,\# B_N]\}$ be the set of
its eigenvalues labeled in increasing order. We shall prove that
\beqn 
\label{trace0}
\Q.\P^\o[X_t=0]=\sup_N \frac 1{\#B_N} \Q\left[\sum_i e^{-\mu_i^\o(B_N)t}\right].
\eeqn 
Let 
\beqnn 
\ce^{\o,N}(f,g)=\frac 12\sum_{\stack{x\sim y}{x,y\in B_N}}
\o(x,y)\, [f(x)-f(y)]\,[g(x)-g(y)]
+\sum_{x\in
B_N}f(x)g(x)\sum_{\stack{y\sim x}{y\notin B_N}}\o(x,y)  
\eeqnn 
be the Dirichlet form of $-{\cal G}^{\o,N}$. From the min-max 
caracterization of $\mu_i^\o(B_N)$, we have  
\beqnn 
\mu_i^\o(B_N)= \max_{f_1,...,f_{i-1}}\min_f \frac {\ce^{\o,N}(f,f)} {\sum_{x\in B_N} f^2(x)}, 
\eeqnn 
where the 'max' is computed on choices of $i-1$ functions defined on $B_N$ and the 'min' is computed 
on functions $f$ such that, for all $j\in\{1,\ldots,i-1\}$, $\sum_{x\in B_N} f(x)f_j(y)=0$. 
For any function $f$, $\ce^{\o,N}(f,f)$ is clearly increasing in $\o$, therefore 
for given $N$, and $i$, $\mu_i^\o(B_N)$ is an increasing function of $\o$ and 
$\sum_i e^{-\mu_i^\o(B_N)t}$ is decreasing in $\o$. Thus (\ref{trace0}) implies the lemma. $\Box$

\vspace{.5cm}
 
\noindent{\bf Proof of (\ref{trace0})} 
Let $\tau_N$ be the exit time of $X_t$ outside $B_N$. 
Note that $\sum_i e^{-\mu_i^\o(B_N)t}$ is just the trace of  
the semi-group of the process $X_t$ killed when leaving the box $B_{N}$, i.e., 
\beqnn 
\sum_i e^{-\mu_i^\o(B_N)t}=\sum_{x\in B_N}\P^\o_x[X_t=x;t<\tau_N].
\eeqnn

We compute $\Q.\P^\o_0[X_t=0]$ using that, from the translation invariance of the probability $\Q$,  
we know that  $\Q.\P^\o_x[X_t=x]$ does not depend on $x$. Therefore  
\begin{eqnarray*}
\Q.\P^\o_0[X_t=0] &=& 
\frac 1{\#B_N} \sum_{x\in B_N} \Q.\P^\o_x[X_t=x]\geq 
\frac 1{\#B_N} \sum_{x\in B_N} \Q.\P^\o_x[X_t=x;t<\tau_N]\\ 
&=& \frac 1{\#B_N} \Q\left[\sum_i e^{-\mu_i^\o(B_N)t}\right] 
\end{eqnarray*} 
proves the lower bound. 

As far as the upper bound is now concerned, note that 
\begin{eqnarray*}
\Q.\P^\o_0[X_t=0] &=& 
\frac 1{\#B_N} \sum_{x\in B_N} \Q.\P^\o_x[X_t=x]\\
&=& \frac 1{\#B_N} \sum_{x\in B_N} \Q.\P^\o_x[X_t=x;t<\tau_{N+k}] 
    + \frac 1{\#B_N} \sum_{x\in B_N} \Q.\P^\o_x[X_t=x;t\geq \tau_{N+k}]\\   
&\leq&  \frac 1{\#B_N} \sum_{x\in B_{N+k}} \Q.\P^\o_x[X_t=x;t<\tau_{N+k}] 
       + \frac 1{\#B_N} \sum_{x\in B_N} \Q.\P^\o_x[t\geq \tau_{N+k}].  
\end{eqnarray*} 

We have 
\beq\nn
\sum_{x\in B_{N+k}} \Q.\P^\o_x[X_t=x;t<\tau_{N+k}] 
= \Q\left[\sum_i e^{-\mu_i^\o(B_{N+k})t}\right]
\leq \#B_{N+k}\, \sup_M \frac 1{\#B_M}\,\Q\left[\sum_i e^{-\mu_i^\o(B_M)t}\right]. 
\eeq

Let $n_t$ be the number of jumps the process $X_t$ performs by time $t$. 
For $x\in B_N$, under $\P^\o_x$, $t\geq \tau_{N+k}$ implies that 
$n_t\geq k$. Therefore 
\begin{eqnarray*}
\sum_{x\in B_N} \Q.\P^\o_x[t\geq \tau_{N+k}] 
&\leq& \sum_{x\in B_N} \Q.\P^\o_x[n_t\geq k]\\ 
&=& \#B_N\, \Q.\P^\o[n_t\geq k], 
\end{eqnarray*} 
using the translation invariance in the last equality. 

So far, we have obtained the bound 
\beqnn
\Q.\P^\o_0[X_t=0] 
\leq \frac {\#B_{N+k}} {\#B_N}\sup_M \frac 1{\#B_M}\Q\left[\sum_i e^{-\mu_i^\o(B_M)t}\right] 
+ \Q.\P^\o[n_t\geq k].
\eeqnn 
First let $N$ tend to $+\infty$, then let $k$ tend to $+\infty$ to deduce that 
\beqnn 
\Q.\P^\o_0[X_t=0] 
\leq \sup_M \frac 1{\#B_M}\Q\left[\sum_i e^{-\mu_i^\o(B_M)t}\right] 
+ \Q.\P^\o[n_t=+\infty].
\eeqnn 

Now the conservativeness assumption and the fact that there are no
instantaneous points of $X_t$ in $\Z^d$ imply that 
$\P^\o[n_t=+\infty]=0$ $\Q$-a.s.  $\Box$


\vskip 1cm
\section{Times of convergence to equilibrium of random walks on the torus}
\setcounter{equation}{0}
\label{sec:times}


Let $S_N$ be the discrete, $d$-dimensional torus of side length 
$N$. When
convenient, we consider $S_N$ as a subset of $\Z^d$. We construct a family of
Markov chains taking their values in $S_N$. Let $\o:S_N\to\R^*_+$, and define 
the Markov generator
\begin{equation}
\label{gen}
\L^{\o,N} f(x)=\sum_{y\sim x}[\o(x)\wedge\o(y)]\,[f(y)-f(x)],
\end{equation} 
where the sum is over sites $y$ which are nearest neighbors to $x$
(relation that is denoted $y\sim x$).
Let $\{X_t,\,t\in\R_+\}$ be the process with distribution $\P_x^{\o,N}$
generated by~(\ref{gen}) and the condition $X_0=x$.
Since $\o(x)>0$ for all $x$, $X_t$ is ergodic under $\P_x^{\o,N}$ for all $x$.
The unique invariant probability measure is the uniform law, denoted by
$\eta_N$. Furthermore, $\eta_N$ is reversible.

We choose the family $\{\o(x),\,x\in\Z^d\}$ i.i.d.~according to a law $\Q$ on
$(\R^*_+)^{\Z^d}$ such that 
\beqn
\label{eq:cond}
\o(x)&\!\!\!\leq\!\!\!&1 \mbox{ for all }x;\\
\Q(\o(0)\leq a)&\!\!\!\sim\!\!\!& a^\g \mbox{ as }a\downarrow0,
\eeqn
where $\g>0$ is a parameter.

\brm
We note that this generator has the same form as 
${\cal G}^\o$ in~(\ref{gen0}) by making $\o(x,y)=\o(x)\wedge\o(y)$,
but for a process in finite volume. We could have defined $\o$ on 
edges, instead of points, as in the previous section,
with i.i.d.~values for different edges, and the same technique 
would apply, with similar results, and heavier computation.
\erm

\brm

If $\o(0)$ were a Bernoulli random variable, then we would have a random walk
on a (independent, site) percolation cluster (provided we started in an
occupied site). See~\cite{kn:MR1}.

\erm

Our main results refer to the following convergence times.
For $\e\in(0,1)$, let
\beqn
\label{t1}
T_1^{\o,N}&=&\inf\{t\geq0:\sup_{x\in S_N}\sup_{|f|\leq1}|\E^{\o,N}_x[f(X_t)]-\eta_N(f)|\leq\e\}\\
\label{t2}
T_2^{\o,N}&=&\inf\{t\geq0:\sup_{|f|\leq1}\sup_{|g|\leq1}|\E^{\o,N}_{\eta_N}[f(X_0)g(X_t)]-\eta_N(f)\eta_N(g)|\leq\e\},
\eeqn
where $\E_x^{\o,N}$ is the expectation with respect to $\P_x^{\o,N}$ and $\E^{\o,N}_{\eta_N}(\cdot)=\int\E_x^{\o,N}(\cdot)\,d\eta_N(x)$.
\brm
The first convergence time is a worst-case one, that is, it is the longest convergence time among
all initial conditions. The second one is an average convergence time among all initial conditions
(under uniform weighting).
\erm
\brm
\label{dom}
Clearly, $T_2^{\o,N}\leq T_1^{\o,N}$ for all $\o$.
\erm

\begin{theo}
\label{t:t1}
For all $\g>0$ and $\e\in(0,1)$, we have $\Q$-a.s.
\beqn
\label{upt1}
\limsup_{N\to\infty}\frac{\log T_1^{\o,N}}{\log N}&\leq&2\vee\frac{d}{\g},\\
\label{lot1}
\liminf_{N\to\infty}\frac{\log T_1^{\o,N}}{\log N}&\geq&2\vee\frac{d}{\g}.
\eeqn
\end{theo}

\begin{theo}
\label{t:t2}
For all $\g>0$ and $\e\in(0,1/4)$, we have $\Q$-a.s.
\beqn
\label{upt2}
\limsup_{N\to\infty}\frac{\log T_2^{\o,N}}{\log N}&\leq&2,\\
\label{lot2}
\liminf_{N\to\infty}\frac{\log T_2^{\o,N}}{\log N}&\geq&2.
\eeqn
In fact, for all $\e\in(0,1/4)$, there exists a constant $c>0$ such that
for all $\o$
\begin{equation}
\label{lo2t2}
\liminf_{N\to\infty}N^{-2}T_2^{\o,N}\geq c.
\end{equation}
\end{theo}

\brm
\label{higho}
If $\o(0)$ were bounded away from zero, that is, if 
$\o(0)>c_1$ $\Q$-a.s.~for some constant $c_1>0$, then 
$\limsup_{N\to\infty}N^{-2}T_1^{\o,N}\leq c_2$ 
$\Q$-a.s.~for some constant $c_2>0$. 
\erm

\brm
Theorems~\ref{t:t1} and~\ref{t:t2} establish that $\Q$-a.s.
\begin{equation}
\label{lims}
\lim_{N\to\infty}\frac{\log T_1^{\o,N}}{\log N}=2\vee\frac{d}{\g}\mbox{ and }
\lim_{N\to\infty}\frac{\log T_2^{\o,N}}{\log N}=2.
\end{equation}
We thus have distinct asymptotic behaviors of $T_1^{\o,N}$ and $T_2^{\o,N}$ 
when $d/\g>2$. 
A heuristic argument to justify that
follows. When $d/\g>2$, $T_1^{\o,N}$, as a worst case convergence time, 
is greater than
or equal to the convergence time starting at a site with minimal $\o$, 
whose order
is clearly smaller than or equal to $N^{-d/\g}$. On the other hand, choosing 
a site 
uniformly at random as a starting point will miss the low $\o$ sites and, 
starting
at high $\o$, the walk will get to equilibrium faster than it will get to any
low $\o$ site. It will be as if there were no low $\o$ sites, and that means
$T_2^{\o,N}$ is of order $N^2$ (see Remark~\ref{higho}).
\erm 

{}From now on, we shall drop the '$N$' in some of our notation. 
For example, we  use the short hand notation $S=S_N$. 


\subsection{Proof of~(\ref{lo2t2})}

Let $A=\{x=(x_1,\ldots,x_d)\in S:x_1\in[0,N/2]\}$, 
$T_A=\inf\{t\geq0:X_t\in A\}$
and, for $\l\geq0$, $h_x^\o(\l)=\E_x^{\o,N}(e^{-\l T_A})$. 
Choosing $f=1_A$ and $g=1_{A^c}$, we have
\beqn
\nonumber
&&\sup_{|f|,|g|\leq1}|\E^{\o,N}_{\eta_N}[f(X_0)g(X_t)]-\eta_N(f)\eta_N(g)|
\geq\eta_N(A)\eta_N(A^c)-\P_{\eta_N}^{\o,N}(X_0\notin A,X_t\in A)\\
\label{lb21}
&\geq&\eta_N(A)\eta_N(A^c)-\P_{\eta_N}^{\o,N}(X_0\notin A,T_A\leq t)
\geq\eta_N(A)\eta_N(A^c)-\inf_{\l>0}\eta_N(1_{A^c}h^\o(\l))e^{\l t}.
\eeqn
We now estimate $\eta_N(1_{A^c}h^\o(\l))$. We will compare with the case
$\o\equiv1$, which corresponds to the usual random walk
on $S_N$.  The Dirichlet form of $X_t$ is given by
\begin{equation}
\label{dir}
\ce^{\o,N}(f,f)=\frac1{2N^d}\sum_{x\sim y\in S}(\o(x)\wedge\o(y))(f(x)-f(y))^2.
\end{equation}
It is clear that $\ce^{\o,N}(f,f)$ is nondecreasing in (the natural partial ordering of) $\o$.
We have also that, for $\l\geq0$,
\begin{equation}
\label{lb22}
\l\eta_N(h^\o(\l))=\inf_{f|_A=1}\ce^{\o,N}(f,f)+\l\eta_N(f^2).
\end{equation}
Since $\ce^{\o,N}(f,f)\leq\ce^\1(f,f)$, where $\1$ is the identically $1$ vector indexed by $S$,
we have that 
\begin{equation}
\label{lb24}
\eta_N(h^\o(\l))\leq\eta_N(h^\1(\l)). 
\end{equation}
Since $T_A$ is a hitting time for an
ordinary rate 1 random walk on $\Z$ under $\P^\1_\cdot$, 
the invariance principle yields
that for all $\l>0$
\begin{equation}
\label{lb23}
\eta_N(h^\1(N^{-2}\l))\to \frac12 +\phi(\l)
\end{equation}
as $N\to\infty$, where $\phi(\l)\to0$ as $\l\to\infty$.
We also have that $\eta_N(h^\o(\l))=\eta_N(A)+\eta_N(1_{A^c}h^\o(\l))$ and $\eta_N(A)\to 1/2$
when $N\to\infty$. Thus, from~(\ref{lb22}),~(\ref{lb24}) and~(\ref{lb23}), 
\begin{equation}
\label{lb25}
\limsup_{N\to\infty}\eta_N(h^\o(N^{-2}\l))
=\frac12+\limsup_{N\to\infty}\eta_N(1_{A^c}h^\o(N^{-2}\l))
\leq\frac12+\phi(\l)
\end{equation}
and it follows that 
\begin{equation}
\label{lb26}
\eta_N(1_{A^c}h^\o(N^{-2}\l))\leq\phi(\l).
\end{equation}
We conclude that
\begin{equation}
\label{lb27}
\liminf_{N\to\infty}\sup_{|f|,|g|\leq1}|\E^{\o,N}_{\eta_N}[f(X_0)g(X_{cN^2})]-\eta_N(f)\eta_N(g)|
\geq\frac14-e\phi(1/c).
\end{equation}
Since $\phi(1/c)\to0$ as $c\to0$, we get that for all $\e<1/4$,
$\liminf_{N\to\infty}N^{-2}T_2^{\o,N}\geq c^*$, where $c^*$ 
is any positive constant
satisfying $\phi(1/c^*)<(1/4-\e)/e$. $\Box$


\subsection{Proof of~(\ref{upt2})}

We make use of {\em generalized Poincar\'e inequalities}~\cite{kn:M}, 
which we recall now.
Let $\B$ denote the set of nearest neighbor bonds of $S$, i.e.,
$\B=\{(x,y):x,y\in S,x\sim y\}$.
For $x,y\in S$, define $r^\o(b)=N^{-d}(\o(x)\wedge\o(y))$, if
$b\in\B$, and $r^\o(b)=0$, otherwise. The Dirichlet form of $\L^{\o,N}$
on $L_2(S,\eta_N)$ can be written as
\begin{equation}
\nn
\ce^{\o,N}(f,f)=\frac1{2}\sum_{b\in\B}(d_bf)^2r^\o(b),
\end{equation}
where $d_bf=f(x)-f(y)$ and the sum ranges over $b=(x,y)$, $x,y\in S$. 

For $p\in(0,2)$, let $q$ be such that $1+1/q=2/p$ and
\begin{equation}
\label{up21}
\tau^{\o,N}(p)=\sup_{f\not\equiv0,\eta_N(f)=0}\,\,
\frac{\eta_N(f^2)^{2/p}}{\ce^{\o,N}(f,f)||f||^{2/q}_\infty}.
\end{equation}
We then have
\begin{equation}
\label{gpi}
T_2^{\o,N}\leq q\e^{-1/q}\tau^{\o,N}(p)
\end{equation}
for all $p\in(0,2)$.

\brm
In the notation of~\cite{kn:M},
$\tau^{\o,N}(p)$, as defined in~(\ref{up21}), equals $1/\k^\o(p)$.
\erm

For all $x,y\in S$, let $\pi_{x,y}$ be a nearest neighbor path from $x$ to $y$
and let $\ell^*=\sup_{x,y}|\pi_{x,y}|$ be the length of the longest path.

Consider now a partitioning of $S=B\cup G$ and let
\begin{equation}
\nn
\tau^{\o,N}_G=\sup_{f\not\equiv0}
\frac{\sum_{b=(x,y)\in G\times G}(d_bf)^2\eta_N(x)\eta_N(y)}
     {\sum_{b\in G\times G}(d_bf)^2r^\o(b)}.
\end{equation}

\begin{lm}
\label{uptau}
\begin{equation}
\nn
\tau^{\o,N}(p)\leq2^{2/q}3^{2/p}\eta_N(B)^{2/p}\ell^*
\sup_{b=(x,y):x\sim y}\frac1{r^\o(b)}+2^{2/q}\tau^{\o,N}_G
\end{equation}
\end{lm}

In the next lemmas, for given $0<\xi<1$, we choose $G$ as the largest 
connected component of the set $\{x:\o(x)\geq\xi\}$ (following a 
deterministic order in case of ties).


\begin{lm}
\label{lowtaug}
For $\xi>0$ small enough, there exists a positive number $c_1$
that depends only on $d$ such that $\Q$-a.s.
\begin{equation}
\nn
\liminf_{N\to\infty}N^2/\tau^{\o,N}_G\geq\xi c_1.
\end{equation}
\end{lm}

\begin{lm}
\label{etab}
For $\xi>0$, there exists a number $c_2(\xi)$
such that $c_2(\xi)\to0$ as $\xi\to0$ and $\Q$-a.s.
\begin{equation}
\nn
\limsup_{N\to\infty}\eta_N(B)\leq c_2(\xi),
\end{equation}
where $B=S\backslash G$.
\end{lm}

\begin{lm}
\label{minw}
There exists a finite number $c_3$ depending only on $d$ 
such that $\Q$-a.s., for all $N$ large enough
\begin{equation}
\nn
\inf_{x\in S}\o(x)\geq N^{-c_3}.
\end{equation}
\end{lm}
In the proof below, we will see that $c_3$ can be taken as $\frac{d}\gamma+\epsilon$ 
for arbitrary $\epsilon>0$.

We postpone the proofs of the above lemmas until after the 
proof of~(\ref{upt2}).

\vspace{.5cm}

\noindent{\bf Proof of~(\ref{upt2}).}
With $\o$ and $\xi>0$ fixed, we choose $N$ big enough so that 
the conclusions of Lemmas~\ref{uptau},~\ref{lowtaug},~\ref{etab} and~\ref{minw}
hold. Then, using also~(\ref{gpi}),
\begin{equation}
\label{eq:upt2a}
T_2^{\o,N}\leq q\e^{-1/q}\tau^{\o,N}(p)\leq
\frac{q\e}4
\{[12\e^{-1}c_2(\xi)]^{2/p}N^{1+d+c_3}+
(c_1\xi)^{-1}(4\e^{-1})^{2/p}N^2]\}.
\end{equation}
Assuming that $\xi$ is small enough, let $p$ satisfy
\begin{equation}
\nn
\frac2p=\frac{(d-1+c_3)\log N+\log(c_1\xi)}{\log c_2(\xi)^{-1}-\log3}.
\end{equation}
With this choice, the two summands in the expression within braces 
in~(\ref{eq:upt2a}) are equal and thus~(\ref{eq:upt2a}) equals
\begin{equation}
\label{eq:upt2c}
\frac{q\e}2(c_1\xi)^{-1}
\exp\frac{\log (c_1\xi)}{\log(c_2(\xi))^{-1}-\log3}
\exp\left\{\left[2+
\frac{(d-1+c_3)\log(4\e^{-1})}{\log(c_2(\xi))^{-1}-\log3}
\right]\log N\right\}.
\end{equation}
Combining~(\ref{eq:upt2a}-\ref{eq:upt2c}), we get
\begin{equation}
\nn
\limsup_{N\to\infty}\frac{\log T_2^{\o,N}}{\log N}\leq
2+\frac{(d-1+c_3)\log(4\e^{-1})}{\log(c_2(\xi))^{-1}-\log3}.
\end{equation}
Since this holds for all $\xi>0$ sufficiently small
and $c_2(\xi)\to0$ as $\xi\to0$,
the result follows.
$\square $


\vspace{.5cm}

\noindent{\bf Proof of Lemma~\ref{uptau}.}
This is very similar to the results of part III in~\cite{kn:MP}.
We estimate the three terms in the decomposition
\begin{equation}
\label{eq:uptaua}
\eta_N(f^2)=
\left(\frac12\sum_{x,y\in G}+\sum_{x\in G,y\in B}+\frac12\sum_{x,y\in B}\right)
(f(x)-f(y))^2\eta_N(x)\eta_N(y)=:I+II+III
\end{equation}
in turn.
\beqn
\nonumber
I&\leq&\frac12(2||f||_\infty)^{2-p}\sum_{x,y\in G}(f(x)-f(y))^p\eta_N(x)\eta_N(y)\\
\nonumber
&\stackrel{\mbox{\small H\"older}}\leq&
2^{1-p}||f||^{2-p}_\infty
\left(\sum_{x,y\in G}(f(x)-f(y))\eta_N(x)\eta_N(y)\right)^{p/2}\\
\label{eq:uptaub}
&\leq&2^{1-p}||f||^{2-p}_\infty
\left(2\ce^{\o,N}(f,f)/\tau^{\o,N}_G\right)^{p/2}.
\eeqn

\beqn
\nonumber
II\!\!\!\!&\leq&\!\!\!\!(2||f||_\infty)^{2-p}\sum_{x\in G,y\in B}(f(x)-f(y))^p\eta_N(x)\eta_N(y)\\
\nonumber
&\stackrel{\mbox{\small H\"older}}\leq&\!\!\!
(2||f||_\infty)^{2-p}
\left(\sum_{x\in G,y\in B}(f(x)-f(y))\eta_N(x)\eta_N(y)\right)^{p}(\eta_N(G)\eta_N(B))^{1-p}\\
\nonumber
&\leq&\!\!\!\!
(2||f||_\infty)^{2-p}
\left(\sum_{x\in G,y\in B}|\!\!\!\!\sum_{\mbox{}\quad b\in\pi_{x,y}}\!\!\!d_bf\,|\,\,\eta_N(x)\eta_N(y)\right)^{p}
\eta_N(B)^{1-p}\\
\nonumber
&\leq&\!\!\!\!
(2||f||_\infty)^{2-p}
\left(\sum_{b}|d_bf|
\sum_{\stack{x\in G,y\in B:}{\pi_{x,y}\ni b}}
\eta_N(x)\eta_N(y)\right)^{p}\!\!\!\eta_N(B)^{1-p}\\
\nonumber
&\stackrel{\mbox{\small H\"older}}\leq&\!\!\!\!
(2||f||_\infty)^{2-p}
\left(\!\sum_{b}|d_bf|^2r^\o(b)\!\right)^{p/2}
\!\left(\!\sum_{b}(r^\o(b))^{-1}\!
\left(\!\sum_{\stack{x\in G,y\in B:}{\pi_{x,y}\ni b}}\eta_N(x)\eta_N(y)\!\right)^{\mbox{}\!\!\!\!2}\,\right)^{p/2}
\!\!\!\eta_N(B)^{1-p}\\
\nonumber
&\leq&\!\!\!\!
(2||f||_\infty)^{2-p}
\left(2\ce^{\o,N}(f,f)\right)^{p/2}
\left(\sup_b(r^\o(b))^{-1}\right)^{p/2}\!\!\!\!
\eta_N(B)^{p/2}\times\\
\nn
&&\mbox{}\hspace{2cm}\times
\left(\sum_{b}\sum_{\stack{x\in G,y\in B:}{\pi_{x,y}\ni b}}\eta_N(x)\eta_N(y)\right)^{p/2}
\!\!\!\!\!\eta_N(B)^{1-p}\\
\nonumber
&\leq&\!\!\!\!
2^{1-p/2}||f||^{2-p}_\infty
\left(\ce^{\o,N}(f,f)\sup_b(r^\o(b))^{-1}\right)^{p/2}
\left(\ell^*\eta_N(B)\right)^{p/2}
\eta_N(B)^{1-p/2}\\
\label{eq:uptauc}
&=&\!\!\!\!
2^{1-p/2}||f||^{2-p}_\infty\eta_N(B)
\left(\ce^{\o,N}(f,f)\sup_b(r^\o(b))^{-1}\ell^*\right)^{p/2},
\eeqn
where the last inequality follows from
\[\sum_{b}\sum_{\stack{x\in G,y\in B:}{\pi_{x,y}\ni b}}\eta_N(x)\eta_N(y)=
\sum_{{x\in G,y\in B}}|\pi_{x,y}|\eta_N(x)\eta_N(y)\leq\ell^*\eta_N(B).\]
Similarly,
\begin{equation}
\label{eq:uptaud}
III\leq
2^{1-p/2}||f||^{2-p}_\infty\eta_N(B)
\left(\ce^{\o,N}(f,f)\sup_b(r^\o(b))^{-1}\ell^*\right)^{p/2}.
\end{equation}

We conclude from~(\ref{eq:uptaua}),~(\ref{eq:uptaub}),~(\ref{eq:uptauc})
and~~(\ref{eq:uptaud}) that
\begin{equation}
\nn
\eta_N(f^2)\leq
\left\{3\eta_N(B)
\left(\sup_b(r^\o(b))^{-1}\ell^*\right)^{p/2}+
\left(\tau^{\o,N}_G\right)^{p/2}\right\}
2^{1-p/2}||f||^{2-p}_\infty(\ce^{\o,N}(f,f))^{p/2}.
\end{equation}
Thus,
\beqn
\nonumber
\tau^{\o,N}(p)
&\leq&
\left\{6\,2^{-p/2}\eta_N(B)
\left(\sup_b(r^\o(b))^{-1}\ell^*\right)^{p/2}+
2^{1-p/2}
\left(\tau^{\o,N}_G\right)^{p/2}\right\}^{2/p}\\
\nonumber
&\leq&
2^{\frac4p-1}\left\{(3\eta_N(B))^{2/p}\sup_b(r^\o(b))^{-1}\ell^*
+\tau^{\o,N}_G\right\}\\
\nn
&=&
2^{\frac2q+1}\left\{(3\eta_N(B))^{2/p}\sup_b(r^\o(b))^{-1}\ell^*
+\tau^{\o,N}_G\right\}.\quad\quad\square
\eeqn
\vspace{.5cm}

\noindent{\bf Proof of Lemma~\ref{lowtaug}.} Since $\o(\cdot)\geq\xi$
on $G$, we have
\begin{equation}
\nn
\tau^{\o,N}_G\leq\frac{\#G}{N^d}\xi^{-1}\tau^{\mathbf\1}_G\leq\xi^{-1}\tau^{\mathbf\1}_G,
\end{equation}
where
\begin{equation}
\nn
\tau^\1_G:=\sup_{f\not\equiv0}
\frac{\sum_{b=(x,y)\in G\times G}(d_bf)^2(\#G)^{-2}}
     {\sum_{b\in G\times G}(d_bf)^2(\#G)^{-1}}
\end{equation}
is the inverse of the spectral gap for the ordinary rate 1 random walk
on $G$.  
{}From Cheeger's inequality, we get that
\begin{displaymath}
\tau^\1_G\leq8\Xi_G^2,
\end{displaymath}
and therefore 
\begin{equation}
\label{eq:lowtaug3}
\tau^{\o,N}_G\leq8\xi^{-1}\Xi_G^2,
\end{equation}
where the isoperimetric constant $\Xi_G$ is defined by: 
\begin{equation}
\nn
\Xi_G:=\sup_{A\subset G}
\frac{\#A\,\#G\setminus A}
     {\#G\,\#\partial_GA},
\end{equation}
where $\partial_GA=\{(x,y):\,x\sim y,x\in A,y\in G\setminus A\}$
is the bond boundary of $A$ with respect to $G$. 
The statement of the Lemma will thus follow if we can prove that 
$N\Xi_G^{-1}$ is bounded from below for large $N$ by some constant that only depends on the
dimension. We shall rather show that 
\begin{equation}
\label{eq:sum}
\sum_N \Q(\Xi_G\geq\a N)<\infty, 
\end{equation} 
for some $\a$. One then uses the Borel-Cantelli Lemma to deduce  from
~(\ref{eq:sum}) that, 
$\Q$.a.s., for large $N$, we have $\Xi_G\leq\a N$ and therefore, 
as follows from~(\ref{eq:lowtaug3}), $\tau^{\o,N}_G\leq8\xi^{-1}\a^2N^2$.  

Following~\cite{kn:MR1}, Subsection 3.1, we note that we can restrict 
ourselves to connected $A$'s such that $G\setminus A$ is connected. 

Since $\#\partial_GA\geq 1$, we have $\frac{\#A\,\#G\setminus
A}{\#G\,\#\partial_GA}\leq\frac\a 2 N$ as soon as $\#A\leq\frac\a 2 N$ 
or $\#G\setminus A \leq\frac\a 2 N$. Thus we may also assume that 
$\#A\geq\frac\a 2 N$ 
and $\#G\setminus A \geq\frac\a 2 N$.

The same argument as in~\cite{kn:MR1}, Subsection 3.1, based on the
classical isoperimetric inequality on $S$, shows that~(\ref{eq:sum}) 
follows from 
\begin{equation}
\label{eq:sum2}
\sum_N \Q\left(\sup_{F\subset\B}\frac{\#F}{\#\{(x,y)\in F; \o(x)\geq\xi,
\o(y)\geq\xi\}}
\geq\a \right)<\infty. 
\end{equation} 
In~(\ref{eq:sum2}), $\B=\{(x,y):x,y\in S,x\sim y\}$ denotes 
the set of nearest neighbor bonds of $S$. The $\sup$ is computed on $*$-connected sets 
$F\subset\B$ such that $\#F\geq \a_1N^{\frac{d-1}d}$, for some constant $\a_1$ that depends on
$\a$ and the dimension.

Given such an $F$, choose a subset, say $\widetilde F$, such that  
$b=(x,y)\ne b'=(x',y')\in\widetilde F\Rightarrow x\ne x'$
and $y\ne y'$. 
Since any point has at most $2d$ neighbors and $\#F\geq \a_1N^{\frac{d-1}{d}}$,
we may assume that $\widetilde F\geq \a_2\#F$, for some positive
$\a_2$.

Now, for all $\l>0$ 
\beqnn
&&\Q(\#\{(x,y)\in F;\o(x)\geq\xi, \o(y)\geq\xi \}\leq\#F/\a)\\ 
&\leq&\Q(\#\{(x,y)\in \widetilde F;\o(x)\geq\xi, \o(y)\geq\xi \}\leq\#F/\a)\\ 
&=& \Q\left(\sum_{(x,y)\in \widetilde F} \1_{\o(x)\geq\xi}\1_{\o(y)\geq\xi}\leq\#F/\a\right)\\ 
&\leq& e^{\frac \l\a \#F} (1-\pi^2+e^{-\l}\pi^2)^{\#\widetilde F}\\
&\leq& e^{\frac \l\a \#F} (1-\pi^2+e^{-\l}\pi^2)^{\a_2\#F}, 
\eeqnn
where $\pi=\Q(\o(x)\geq\xi)$. 

By the above inequality, and the fact that
the number of distinct $*$-connected subsets $F$ with $\#F=n$ 
is bounded above by $N^d e^{\a_3n}$ for some 
$\a_3$~\cite{kn:S}, we get
\beqn 
\nonumber
\Q\left(\sup_F\frac{\#F}{\#\{(x,y)\in F;\o(x)\geq\xi,\o(y)\geq\xi\}}\geq\a\right)
&\leq& N^d\sum_{n\geq \a_1N^{\frac{d-1}{d}}}
e^{[\a_3+\l\a^{-1}+\a_2\log(1-\pi^2+e^{-\l}\pi^2)]n}\\
&=&
N^d\sum_{n\geq \a_1N^{\frac{d-1}{d}}}
e^{-\a_4n},
\nn
\eeqn
where $\a_4:=-[\a_3+\l\a^{-1}+\a_2\log(1-\pi^2+e^{-\l}\pi^2)]>0$,
provided we choose $\l$ and $\a$ such that 
$\a_3+\l/\a<\l\a_2$  
and $\xi\leq\xi_0$, for $\xi_0$ close enough 
to $0$, depending on $\a,\l,\a_2,\a_3$ and $\g$ only. $\square$

\vspace{.5cm}

\noindent{\bf Proof of Lemma~\ref{etab}.}
Consider the site percolation model on $\Z^d$ where a site $x$
is occupied if $\o(x)\geq\xi$. Let $\xi_0$ be positive and 
satisfy $\Q(\o(x)\geq\xi_0)>p_c$, the critical density
for the a.s.~appearance of an infinite connected component $C$.
Then, if $\xi<\xi_0$, $C$ exists a.s.
Let $\tilde C_N=C\cap\tilde S_N$, where $\tilde S_N$ is $S_N$ viewed
as a subset of $\Z^d$ 
(that is, without the {\em boundary identification}), say,
$\tilde S_N=(-N/2,N/2]^d\cap\Z^d$. 
Let $C_N$ be $\tilde C_N$ viewed as a subset of the torus $S_N$
(that is, with the boundary identification).
Then, it follows by standard ergodicity
arguments that $\lim_{N\to\infty}\eta_N(\tilde C_N)=\theta(\xi):=\Q(0\in C)$
$\Q$-a.s.
Since $\theta(\xi)\to1$ as $\xi\to0$ (a well known result~\cite{kn:G}),
the result would follow if $C_N$ were connected, which it is not
necessarily. 

Consider then $\hat C_{N}:=\tilde C_{N-\lfloor\sqrt N\rfloor}$.
We claim that $\hat C_{N}$ is connected in $\tilde S_N$,
and thus also in $S_N$, for all large enough $N$ $\Q$-a.s.
Indeed, in the event that $\hat C_{N}$ is not connected in $\tilde S_N$,
there exist two sites at the boundary of $\tilde S_{N-\lfloor\sqrt N\rfloor}$
that are connected to the boundary of $\tilde S_{N}$ but are not connected
to one another. This implies that there exists a site $\tilde x$ 
at the boundary of $\tilde S_{N}$ whose (occupied) cluster 
(in $\tilde S_{N}$) has a boundary (of vacant sites) of size at least 
$\lfloor\sqrt N\rfloor$.
Now, the (bond) boundary of any finite cluster of a site in $\tilde S_{N}$ 
can be identified with a surface of {\em plaquettes} around the given site,
each plaquette crossing orthogonally a boundary bond. For each
such plaquette, there corresponds thus an inner occupied site and
an outer vacant one. For a given such surface of plaquettes of size
(total number of plaquettes) $n$, there is at least $n/(2d)$ 
distinct outer vacant sites (since a vacant site can not be adjacent
to more than $2d$\footnote{Actually, $n/(2d-1)$ is a better bound.}
In the case of $\tilde x$, the surface of plaquettes will intersect the 
boundary of $\tilde S_{N}$ in a closed curve. It will also have to cross
the region between the boundaries of $\tilde S_{N}$ and 
$\tilde S_{N-\lfloor\sqrt N\rfloor}$. For this reason it will contain 
at least $\lfloor\sqrt N\rfloor$ plaquettes. 

{}From the arguments in the latter paragraph, we get the following estimate.
\beqn
\nonumber
&\Q(\hat C_{N}\mbox{ is not connected in }\tilde S_N)&\\
\label{eq:etab1}
&\leq
\sum_{\tilde x\in\partial\tilde S_N}
\sum_{\stack{\Gamma\,\mbox{{\tiny around}}\,x:}
      {\#\Gamma\geq\lfloor\sqrt N\rfloor}}
\Q(\mbox{all the sites at the outer boundary of $\Gamma$ are vacant}),& 
\eeqn
where the latter sum above is over surface of plaquettes $\Gamma$ around 
$\tilde x$. The number of distinct such surfaces which have size $n$
can be estimated to be exponential in $n$~\cite{kn:S}. Proceeding
with the estimation we get that the
right hand side of~(\ref{eq:etab1}) equals
\beqn
\nonumber
&\sum_{\tilde x\in\partial\tilde S_N}
\sum_{n\geq\lfloor\sqrt N\rfloor}
\sum_{\stack{\Gamma\,\mbox{{\tiny around}}\,x:}
     {|\Gamma|=n}}
\Q(\mbox{all the sites at the outer boundary of $\Gamma$ are vacant})&\\
\nn
&\leq
N^d\sum_{n\geq\lfloor\sqrt N\rfloor}e^{\nu n}[\Q(\o(0)<\xi)]^{n/(2d)},&
\eeqn
where $\nu$ depends only on $d$. Thus, by taking $0<\xi<\xi_0$ small
enough, the probability in the left hand side of~(\ref{eq:etab1})
can be made summable and the claim at the beginning of the previous
paragraph follows by Borel-Cantelli. The lemma then follows. $\square$ 

\vspace{.5cm}

\noindent{\bf Proof of Lemma~\ref{minw}.}
We will prove that $\Q$-a.s.
\beqn
\nn
\lim_{N\to\infty}\frac{\log\inf_x\o(x)}{\log N}=-\frac d\g.
\eeqn
For that, let $c<d/\g$. Then
\beqn
\label{eq:minw2}
\Q(\inf_x\o(x)\geq N^{-c})=[\Q(\o(x)\geq N^{-c})]^{N^d}
\leq(1-c_1N^{-c\g})^{N^d}\leq e^{-c_1N^{d-c\g}}
\eeqn
for $N$ large enough and some constant $c_1$. 
Thus the Borel-Cantelli lemma implies the upper bound in~(\ref{eq:minw2}). 

Now, let $c>d/\g$. For $e^k\leq N\leq e^{k+1}$, we have
\beqn
\nn
\inf_{x\in S_N}\o(x)\geq
\inf_{x\in S_{e^k}}\o(x)\wedge\inf_{x\in S_{e^{k+1}}\setminus S_{e^k}}\o(x).
\eeqn
Therefore,
\beqn
\nonumber
&&
\Q\left(\exists N\in[e^k,e^{k+1}):\,\inf_{x\in S_N}\o(x)\leq N^{-c}\right)\\
\nonumber
&\leq&\Q\left(\inf_{x\in S_{e^k}}\o(x)\leq e^{-ck}\right)+
  \Q\left(\inf_{x\in S_{e^k}\setminus S_{e^k}}\o(x)\leq e^{-ck}\right)\\
\label{eq:minw4}
&=&(1-(1-c_1e^{-c\g k})^{e^{dk}})+(1-(1-c_1e^{-c\g k})^{e^{d(k+1)}-e^{dk}})
\leq c_2\,e^{-(c\g-d)k}
\eeqn
and the result follows from Borel-Cantelli and the summability of the 
probabilities on the left hand sides of~(\ref{eq:minw2}) 
and~(\ref{eq:minw4}), implied by their right hand sides. $\square$


\subsection{Proof of~(\ref{lot1})}

{From} $T_1^{\o,N}\geq T_2^{\o,N}$ and $\liminf_{N\to\infty} N^{-2}T_2^{\o,N}>c$ $\Q$-a.s., 
we deduce that $\liminf_{N\to\infty} N^{-2}T_1^{\o,N}>c$ $\Q$-a.s. and, thus,
$\liminf_{N\to\infty}{\log T_1^{\o,N}}/{\log N}\geq2$ $\Q$-a.s.

We argue now for the inequality  
$\liminf_{N\to\infty}{\log T_1^{\o,N}}/{\log N}\geq d/\g$
$\Q$-a.s. Let $x\in S$. During an exponential time of parameter
$\sum_{y:y\sim x}\o(y)\wedge\o(x)$, the process $X$ starting at $x$ 
stays still. Therefore,
\beqn\nn
\sup_{|f|\leq1}|\E_x^{\o,N} f(X_t)-\eta_N(f)|&\geq&\P_x^{\o,N}(X_t=x)-N^{-d}\\
\nn
&\geq& e^{-t\sum_{y\sim x}\o(y)\wedge\o(x)}-N^{-d}
\geq e^{-2d\o(x)t}-N^{-d},
\eeqn
i.e., 
\beqn
\nn
T_1^{\o,N}\geq\frac1{2d}\sup_x\o(x)^{-1}\log(\e+N^{-d})^{-1}.
\eeqn
Therefore,
\beqn
\nn
\frac{\log T_1^{\o,N}}{\log N}\geq 
\frac{\log\sup_x\o(x)^{-1}}{\log N}+o(1).
\eeqn

Now, let $0<\d<1$ be arbitrary.
\beqn
\nn
\Q\left(\log\sup_x\o(x)^{-1}\leq(1-\delta)\frac d\g\log N\right)=
[\Q(\o(x)\geq N^{-(1-\d)d/\g})]^{N^d}
\leq[1-N^{-(1-\d')d}]^{N^d},
\eeqn
for any $1>\d'>\d$, provided $N$ is large enough. Thus, the above probability
is summable in $N$ for any $\d>0$, and the result follows by Borel-Cantelli.
$\square$


\subsection{Proof of~(\ref{upt1}): Spectral gap estimates}
Let
\begin{equation}
\nn
\tau^{\o,N}=\sup_{f\not\equiv0,\eta_N(f)=0}\,\,
\frac{\eta_N(f^2)}{\ce^{\o,N}(f,f)}
\end{equation}
be the inverse of the spectral gap. From general facts~\cite{kn:SC}, we have
\begin{equation}
\nn
|\E_x^{\o,N}[f(X_t)]-\eta_N(f)|\leq\eta_N(x)^{-1/2}e^{-t/\tau^{\o,N}}, 
\end{equation}
where $f$ is any function uniformly bounded by $1$. 
Thus
\begin{equation}
\nn
\limsup_{N\to\infty}\frac{\log T_1^{\o,N}}{\log N}\leq
\limsup_{N\to\infty}\frac{\log\tau^{\o,N}}{\log N}.
\end{equation}

Using a formula of Saloff-Coste (see Theorem 3.2.3 in~\cite{kn:SC}), we get
\begin{equation}
\label{upt1d}
\tau^{\o,N}\leq N^{-d}\max_{b\in\B}\frac{W(b)}{\o(b)}
\sum_{\stack{(x,y):}{\pi_{x,y}\ni b}}|\pi_{x,y}|_{{}_W}
=N^{-d}\max_{b\in\B}\frac{W(b)}{\o(b)}
\sum_{b'\in\B}\frac1{W(b')}\,\n(b,b'),
\end{equation}
where $\o(b')=\o(x')\wedge\o(y')$ for $b'=(x',y')\in\B$, $W:\B\to(0,\infty)$
is an arbitrary weight function,
$\{\pi_{x,y}:(x,y)\in S\times S\}$ is an arbitrary complete set of paths
($\pi_{x,y}$ is a path with end points $x$ and $y$), for an arbitrary path
$\pi$ in $S$,
$|\pi|_{{}_W}=\sum_{b\in\pi}1/W(b)$,
and $\n(b,b'):=\#\{(x,y)\in S\times S:\,b,b'\in\pi_{x,y}\}$.

It remains to estimate the right hand side of~(\ref{upt1d}).
The key point here is the choices of the weight function and the complete set 
of paths. Roughly speaking, the latter will be taken in such a way that no 
path in it has {\em interior sites} with low values of $\o$; and the former 
will give low weight to bonds with low values of $\o$. We are precise next.

\begin{df}
  \label{df:gp}
  Given $\e>0$, a site $x\in S$ will be called $\e$-{\em good} if
  $\o(x)>N^{-\e}$. Otherwise, it will be called $\e$-{\em bad}.
  A bond $b=(x,y)\in\B$ will be $\e$-{\em good} if $x$ and $y$ 
  are $\e$-{\em good}. Otherwise, it will be called $\e$-{\em bad}.
\end{df}

\begin{df}
  \label{df:ip}
  Given $L>0$ and a path $\pi\in S$ connecting given sites $x,y$, 
  a site
  $z$ in $\pi$ will be called an $L$-{\em interior site} of $\pi$
  if $||z-x||_\infty,||z-y||_\infty>L$.
\end{df}

\begin{df}
  \label{df:gsp}
  Given $L,\e>0$ and $\Gamma$, a set of paths of $S$, $\Gamma$
  will be called $(L,\e)$-good if all the paths of $\Gamma$
  have all their $L$-interior sites, if any, $\e$-good.  
\end{df}

We now construct for every $N$ a complete set of paths for $S_N$
which will turn out to be almost surely $(L,\e)$-good for all 
large enough $N$ and which will have other properties leading to 
the validity of~(\ref{upt1}). 

We start with an auxiliary set of paths.
\begin{df}
\label{aux}
For $x,y\in S$, let $\eta_{x,y}$ be the path given by moving sequentially
in the $1$-st, $2$-nd,..., $d$-th coordinate direction one step at a time,
along the longest segment (and according to an arbitrary predetermined
order in case of a tie), from $x$ to $y$,
until the coordinates are successively matched.
\end{df}
For example, if $d=3$, $N=100$, $S_N=\{1,2,\ldots,100\}^3$ (with the
boundaries appropriately identified), $x=(1,1,1)$ and $y=(2,20,80)$,
then $\eta_{x,y}=\gamma_1\cup\gamma_2\cup\gamma_3$ is the union of the
segments
\begin{eqnarray*}
  \gamma_1&=&\{(1,1,1)\equiv(100,1,1),(99,1,1),\ldots,(3,1,1),(3,1,1)\}\\
  \gamma_2&=&\{(2,1,1)\equiv(2,100,1),(2,99,1),\ldots,(2,21,1),(2,20,1)\}\\
  \gamma_3&=&\{(2,20,1),(2,20,2),\ldots,(2,20,79),(2,20,80)\}.
\end{eqnarray*}
Now for $L>0$ 
we define the $L$-{\em sausage} $\S_L=\S_L(x,y)$ 
with base $\eta_{x,y}$ and width $L$ as follows. We suppose $N>3L$. 
Let $i_1,i_2,\ldots,i_k$, $1\leq k\leq d$ 
be the coordinates where $x$ differs from $y$ in increasing order, so
that $\eta_{x,y}$ is the union of the segments $\gamma_1,\ldots,\gamma_k$,
each of length at least $N/2$, with $\gamma_i$ parallel to the coordinate
direction $i$. 
If $k<d$, then let 
$i^\ast=\min\{i:1\leq i\leq d \mbox{ and }i\ne i_k, 1\leq k\leq d\}$
and 
\begin{eqnarray}
  \nonumber
  \S_L&=&\{(z_1,\ldots,z_{i^\ast-1},w_{i^\ast},z_{i^\ast+1},\ldots,z_d):\\
  \nn
  &&
  \hspace{1cm} 
  z_{i^\ast}\leq w_{i^\ast}\leq z_{i^\ast}+L-1,
  (z_1,\ldots,z_{i^\ast-1},z_{i^\ast},z_{i^\ast+1},\ldots,z_d)\in\eta_{x,y}\}.
\end{eqnarray}
If $k=d$, then let 
\begin{eqnarray}
  \nn
  \S'_L&=&\{(w_1,z_2\ldots,z_d): z_1\leq w_1\leq z_1+L-1,
  (z_1,\ldots,z_d)\in\cup_{j=2}^k\gamma_j\},\\
  \nn
  \S''_L&=&\{(w_1,z_2\ldots,z_d): z_1-L+1\leq w_1\leq z_1,
  (z_1,\ldots,z_d)\in\cup_{j=2}^k\gamma_j\}.
\end{eqnarray}
Now let $R_1$ be the uniquely defined rectangle with base $\gamma_1$
and width $L$ such that either $R_1\cap\S'_L$ or $R_1\cap\S''_L$
is a $L\times L$ square (one and only one of these possibilities occurs).
In the latter case, $\S_L=R_1\cup\S''_L$; in the former one,
$\S_L=R_1\cup\S'_L$.

\begin{figure}[!ht]
\begin{center}
\includegraphics[width=7cm]{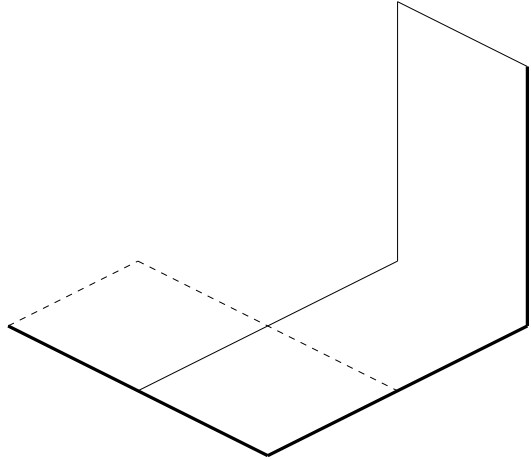}
\caption{Thick polygonal is $\eta_{x,y}$; rectangle delimited by
dashed lines is $R_1$; strip delimited by thin lines
is $\S''_L$}\label{sausage}
\label{sau} 
\end{center}
\end{figure}

\begin{rmk}
\label{rmk:strip}
Notice that $\S_L$ can be seen as either a single 
bidimensional\footnote{Even if living in $k$-dimensional space.} strip 
of length at least $N/2$ and at most $dN$ and width $L$, when $k<d$, or 
the union of two such strips (one of which is the rectangle $R_1$), 
when $k=d$. 
\end{rmk}

Given $\e>0$ and a strip $\S$ of length at least $N/2$ and 
at most $dN$ and width $L$, 
we consider the site percolation model in $\S$ in which a site is open 
if and only if it is $\e$-good and define the event $A_{\S}=A_{\S}(L)$ 
that there 
exists an open path connecting the two smaller sides of $\S$
(within $\S$).
Then one argues as usually that $A^c_{\S}$ 
is the event that there exists a $\ast$-closed path
connecting the two larger sides of $\S$ (within $\S$).
It is clear that $A_{\S}(L)\subseteq A_{\S}(L')$ if $L\leq L'$.

Now consider the event $A_N=A_N(L)$ that $A_{\S}$ occurs for all the strips
involved in the sausages $\S_L(x,y)$ for all $x,y\in S$. 
Clearly, $A_N(L)\subseteq A_N(L')$ if $L\leq L'$.
\begin{df}
\label{ell}
Let
\begin{equation}
  \nn
  \ell_\e=\inf\{L:\,3L<N\mbox{ and }A_N(L)\mbox{ occurs}\},
\end{equation}
with the convention that $\inf\emptyset=\infty$.
\end{df}

The following result will be proven below.

\begin{prop}
  \label{prop:bound}
\begin{equation}
  \label{eq:bound}
  \tau^{\o,N}\leq C(\ell_\e+1)^{2d}\left(N^{2+\e}+\max_{x\in S}\frac1{\o(x)}\right),
\end{equation}
where $C>0$ depends only on $d$.
\end{prop}

This (deterministic) result, together with the following (probabilistic) one 
yields~(\ref{upt1}), after one uses Lemma~\ref{minw} and Borel-Cantelli.

\begin{lm}
  \label{lm:prob}
  For all large enough $N$
\begin{equation}
  \nn
  \P\left(\ell_\e>\left\lceil 4\frac{d+1}{\g\e}\right\rceil\right)
  \leq\frac c{N^{1+\d}},
\end{equation}
where $c$ depends only on $d$ and $\d>0$ is independent of $N$.
\end{lm}

\noindent{\bf Proof of Lemma~\ref{lm:prob} }
For $L>0$ fixed, we have that
\begin{equation}
  \nn
  \P(A_{\S}^c(L))\leq 
  dN\max_{x\in\bar\S}\P(x\mbox{ is connected within $\S$ by 
  a $\ast$-closed path to }\underline\S),
\end{equation}
where $\bar\S$ and $\underline\S$ are the two larger sides of $\S$.
Now the latter probability can be bounded above in a standard way by
\begin{equation}
  \nn
  \sum_{l\geq L}\lambda_l N^{-\g\e l},
\end{equation}
where $\lambda_l$ is the number of distinct $\ast$-paths of length $l$
within $\S$ and starting at $x$. This is bounded above in a standard way 
by $7^l$ and thus
\begin{equation}
  \nn
  \P(A_{\S}^c(L))\leq dN\sum_{l\geq L}(7N^{-\g\e})^{l}\leq c N^{1-\g\e L/2},
\end{equation}
for some constant $c$ and all large enough $N$.

Then 
\begin{equation}
  \nn
  \P(A_N^c(L))\leq  c N^{2d+1-\g\e L/2}.
\end{equation}

The result now follows from the observation that 
$\{\ell_\e>L\}\subset A_N^c(L)$. $\Box$

\subsection{Proof of Proposition~\ref{prop:bound}}

We assume $\ell_\e<\infty$; otherwise, the bound is obvious.
We choose the weight function $W$.
For $b\in\B$, we make
\begin{equation}
\nn
W(b)=
\begin{cases}
1,& \mbox{if } b \mbox{ is $\e$-good},\\
\frac1N,& \mbox{if } b \mbox{ is $\e$-bad}. 
\end{cases}
\end{equation}

We now choose a complete set of paths for $S$, $\Gamma$.
Since $\ell_\e<\infty$,
we have that for all $x,y\in S$, there will be a 
$(\leps,\e)$-good path within $\S_{\leps}(x,y)$ connecting $x$ and $y$, 
so we choose one of them 
(according to some arbitrary predetermined order), 
call it $\pi_{x,y}$, and make
\begin{equation}
  \nn
  \Gamma=\{\pi_{x,y}; x,y\in S\}.
\end{equation}

We now use the above $W$ and $\Gamma$ in~(\ref{upt1d}).
Let $\B_1=\{b\in\B: b$ is $\e$-good\} and
$\B_2=\{b\in\B: b$ is $\e$-bad\} $=\B\setminus\B_1$.  Then
\begin{equation}
\label{sc1}
\tau^{\o,N}\leq\tau^\o_{11}+\tau^\o_{12}+\tau^\o_{21}+\tau^\o_{22},
\end{equation}
where, for $i,j=1,2$,
\begin{equation}
\nn
\tau^\o_{ij}=N^{-d}\max_{b\in\B_i}\frac{W(b)}{\o(b)}
\sum_{b'\in\B_j}\frac1{W(b)}\,\n(b,b').
\end{equation}

For $x,y\in S$, let $Q_x$, resp.~$Q_y$, denote the $\leps\times\leps$ square 
contained in $\S_{\leps}(x,y)$ with $x$, resp.~$y$, as one of its corners.

\begin{rmk}
\label{good}
Notice that for every $x,y\in S$, all bonds of 
$\pi_{x,y}\setminus(Q_x\cup Q_y)$ are $\e$-good. 
\end{rmk}

Given $b,b'\in\B$, let 
$\m(b,b')=\#\{(x,y)\in S\times S:\,b,b'\in\eta_{x,y}\}$; see
Definition~\ref{aux}.

\medskip

\noindent{\bf Estimation of $\tau^\o_{11}$.} 
\begin{equation}
\nn
\tau^\o_{11}\leq N^{\e-d}\max_{b\in\B}
\sum_{b'\in\B}\n(b,b').
\end{equation}
Now for every $b,b'\in\B$
\begin{eqnarray}
\nn
\n(b,b')&\leq&\#\{(x,y)\in S\times S:\,b,b'\in\S_L(x,y)\}\\
\nn
&\leq&\#\{(x,y)\in S\times S:\,a,a'\in\eta_{x,y}
\mbox{ for some }a,a'\in\B:
\mbox{dist}(a,b)\vee\mbox{dist}(a',b')\leq \leps\}\\
\nn
&\leq&\sum_{a,a':\,
\mbox{\scriptsize dist}(a,b)\vee\mbox{\scriptsize dist}(a',b')\leq\leps}
\m(a,a').
\end{eqnarray}
where dist is the usual Hausdorff distance between sets. Thus
\begin{equation}
\label{sc5}
\tau^\o_{11}\leq\leps^{2d}N^{\e}M_N
\end{equation}
where $M_N:=N^{-d}\max_{a\in\B}\sum_{a'\in\B}\m(a,a')$.

To estimate $M_N$, we start with the observation that 
since our paths are described in an oriented way, we must specify
which of $a$ or $a'$ is traversed first and in which direction.
Given $a=(w,z)$, we have
\begin{eqnarray}
\label{mn1a}
\sum_{a'=(w',z')\in\B}\m(a,a')
&=&\sum_{a'\in\B}\#\{(x,y):\,a,a'\in\eta_{x,y}
\mbox{ in the order }w,z,w',z'\}\\
\label{mn1b}&+&\sum_{a'\in\B}\#\{(x,y):\,a,a'\in\eta_{x,y}
\mbox{ in the order }z,w,w',z'\}\\
\label{mn1c}&+&\sum_{a'\in\B}\#\{(x,y):\,a,a'\in\eta_{x,y}
\mbox{ in the order }w',z',w,z\}\\
\label{mn1d}&+&\sum_{a'\in\B}\#\{(x,y):\,a,a'\in\eta_{x,y}
\mbox{ in the order }w',z,'z,w\}.
\end{eqnarray}

We estimate the sum in~(\ref{mn1a}). The estimation for the ones 
in~(\ref{mn1b}-\ref{mn1d}) is similar.
Let $j$ be the coordinate where $w,z$ differ, that is
$z_i=w_i$ if $i\ne j$ and $z_j=w_j\pm 1$.
Then the ordering imposes that $z_i'=w_i'=w_i$ if $i<j$.
The sum in~(\ref{mn1a}) 
can then be decomposed as follows.
\begin{equation}
\label{mn2}
\sum_{k=j}^d\sum_{a'\in\Lambda_k}\m'(a,a'),
\end{equation}
where $\m'(a,a')=\#\{(x,y):\,a,a'\in\eta_{x,y}\mbox{ in the order }w,z,w',z'\}$
and 
\begin{equation}
\nn
\Lambda_k=\{(w',z')\in\B: z_i'=w_i'=w_i,\mbox{ if } i<j;
 z_i'\ne w_i,\mbox{ if } j\leq i\leq k;
z_i'=w_i,\mbox{ if }k<i\leq d\}. 
\end{equation}
It is clear that $|\Lambda_k|\leq N^{k-j+1}$. Now, for $a'\in\Lambda_k$
\begin{eqnarray}
\nn
\m'(a,a')&\leq&\#\{x\in S:x_i=w_i\mbox{ for } i>j\}
\times\#\{y\in S:y_i=z_i'\mbox{ for } i<k\}\\
\nn
&\leq& N^j N^{d-k+1}.
\end{eqnarray}
Thus~(\ref{mn2}) and~(\ref{mn1a}) are bounded above by
$dN^{2+d}$. After a similar reasoning for~(\ref{mn1b}-\ref{mn1d}),
with the same bounds, we finally get from~(\ref{sc5}) that
\begin{equation}
\label{sc6}
\tau^\o_{11}\leq 4d\leps^{2d}N^{2+\e}.
\end{equation}

\medskip

\noindent{\bf Estimation of $\tau^\o_{12}$.} 
\begin{equation}
\nn
\tau^\o_{12}\leq N^{\e-d+1}\max_{b\in\B}
\sum_{b'\in\B_2}\n(b,b').
\end{equation}
By Remark~\ref{good}, if $b'\in\B_2$ is in $\pi_{x,y}\in\Gamma$, then
$b'$ must be either in $Q_x$ or in $Q_y$ (see definition right
above Remark~\ref{good}). Thus
\begin{eqnarray}
\nn
\n(b,b')&\leq&\#\{(x,y)\in S\times S:\,a\in\eta_{x,y}
\mbox{ for some }a\in\B\mbox{ and }
\mbox{dist}(a,b)\vee\mbox{dist}(x,b')\leq\leps\}\\
\nn
&+&\#\{(x,y)\in S\times S:\,a\in\eta_{x,y}
\mbox{ for some }a\in\B\mbox{ and }
\mbox{dist}(a,b)\vee\mbox{dist}(y,b')\leq\leps\}\\
\nn
&\leq&\sum_{a\in\B,z\in S:\,
\mbox{\scriptsize dist}(a,b)\vee\mbox{\scriptsize dist}(z,b')\leq\leps}
[\j(a,z)+\tilde\j(a,z)],
\end{eqnarray}
where
\begin{equation}
\nn
\j(a,z)=\#\{x\in S:\,a\in\eta_{x,z}\},\quad
\tilde\j(a,z)=\#\{y\in S:\,a\in\eta_{z,y}\}.
\end{equation}
We conclude that
\begin{eqnarray}
\nn
\tau^\o_{12}
&\leq&\leps^{d}N^{\e-d+1}\left[\max_{a\in\B}\sum_{z\in S}\j(a,z)\w(z)
+\max_{a\in\B}\sum_{z\in S}\tilde\j(a,z)\w(z)\right]\\
\label{sc9a}
&\leq&\mbox{const }\leps^{2d}N^{\e-d+1}
\left[\max_{a\in\B}\sum_{z\in S}\j(a,z)
+\max_{a\in\B}\sum_{z\in S}\tilde\j(a,z)\right],
\end{eqnarray}
since
\begin{equation}
\nn
\w(z):=\#\{b'\in\B:\mbox{dist}(z,b')\leq\leps\}\leq\mbox{const }\leps^{d}.
\end{equation}

We estimate the first max term in~(\ref{sc9a}). 
The other one is treated
similarly, with the same bound. Let $a=(u,v)$. We decompose
$\j(a,z)$ in $\j'(a,z)$ and $\j''(a,z)$, where
\begin{eqnarray}
\nn
\j'(a,z)&=&
\#\{x\in S:\,a\in\eta_{x,z},\mbox{ with $u$ traversed before $v$}\},\\
\nn
\j''(a,z)&=&
\#\{x\in S:\,a\in\eta_{x,z},\mbox{ with $v$ traversed before $u$}\}.
\end{eqnarray}
We estimate $\max_{a\in\B}\sum_{z\in S}\j'(a,z)$. 
The expression involving $\j''(a,z)$
is treated similarly, with the same bound.
Let $j$ be the coordinate where $u$ and $v$ differ. Then $z$ must 
satisfy $z_i=u_i$, if $1\leq i\leq j-1$. We conclude that there
are at most $N^{d-j+1}$ such $z$'s. For each one, if $a\in\eta_{x,z}$,
then $x$ must satisfy $x_i=u_i$, if $j+1\leq i\leq d$.
We conclude that there are at most $N^{j}$ such $x$'s.
Thus,  
\begin{equation}
\nn
\max_{a\in\B}\sum_{z\in S}\j'(a,z)
\leq\max_{1\leq j\leq d}N^{j}\,N^{d-j+1}=N^{d+1}.
\end{equation}
We conclude that
\begin{equation}
\label{sc13}
\tau^\o_{12}\leq\mbox{const }\leps^{2d}N^{2+\e}.
\end{equation}

\noindent{\bf Estimation of $\tau^\o_{21}$.} 
\begin{equation}
\nn
\tau^\o_{21}\leq\left(\max_{x\in S}\frac1{\o(x)}\right)
                N^{-d-1}\max_{b\in\B_2}\sum_{b'\in\B}\n(b,b').
\end{equation}

We now estimate the max of the sum above, in much the same way 
as we estimated $\max_{b\in\B}\sum_{b'\in\B_2}\n(b,b')$ above.
By Remark~\ref{good}, if $b\in\B_2$ is in $\pi_{x,y}\in\Gamma$, then
$b$ must be either in $Q_x$ or in $Q_y$. Thus
\begin{eqnarray}
\nn
\n(b,b')&\leq&\#\{(x,y)\in S\times S:\,a\in\eta_{x,y}
\mbox{ for some }a\in\B\mbox{ and }
\mbox{dist}(a,b')\vee\mbox{dist}(x,b)\leq\leps\}\\
\nn
&+&\#\{(x,y)\in S\times S:\,a\in\eta_{x,y}
\mbox{ for some }a\in\B\mbox{ and }
\mbox{dist}(a,b')\vee\mbox{dist}(y,b)\leq\leps\}\\
\nn
&\leq&\sum_{a\in\B,z\in S:\,
\mbox{\scriptsize dist}(a,b')\vee\mbox{\scriptsize dist}(z,b)\leq\leps}
[\j(a,z)+\tilde\j(a,z)].
\end{eqnarray}
Thus,
\begin{eqnarray}
\nn
\max_{b\in\B_2}\sum_{b'\in\B}\n(b,b')
&\leq&\mbox{const }\leps^{d}\left[\max_{z\in S}\sum_{a\in\B}\j(a,z)\bar\w(a)
+\max_{z\in S}\sum_{a\in\B}\tilde\j(a,z)\bar\w(a)\right]\\
\label{sc16a}
&\leq&\mbox{const }\leps^{2d}\left[\max_{z\in S}\sum_{a\in\B}\j(a,z)
+\max_{z\in S}\sum_{a\in\B}\tilde\j(a,z)\right],
\end{eqnarray}
where 
\begin{equation}
\nn
\bar\w(a):=\#\{b'\in\B:\mbox{dist}(a,b')\leq\leps\}\leq\mbox{const }\leps^{d}.
\end{equation}
We estimate the first summand within square brackets in~(\ref{sc16a}). The
second one can be similarly estimated with the same resulting bound.
\begin{equation}
\label{sc18}
\max_{z\in S}\sum_{a\in\B}\j(a,z)\leq\max_{z\in S}\sum_{a\in\B}\j'(a,z)
+\max_{z\in S}\sum_{a\in\B}\tilde\j''(a,z)
\end{equation}
and we estimate the first summand within square brackets in~(\ref{sc18}) only.
The second one can be similarly treated with the same bound.
Let $z\in S$ be fixed and $j$ be the coordinate where $u$ and $v$ differ,
where $(u,v)=a$. Then $u$ must 
satisfy $u_i=z_i$, if $1\leq i\leq j-1$. We conclude that there
are at most $N^{d-j+1}$ such $u$'s. For each one, if $a\in\eta_{x,z}$,
then $x$ must satisfy $x_i=u_i$, if $j+1\leq i\leq d$.
We conclude that there are at most $N^{j}$ such $x$'s.
We then conclude that
\begin{equation}
\nn
\max_{z\in S}\sum_{a\in\B}\j'(a,z)
\leq \max_{z\in S}\sum_{j=1}^d
\sum_{\stack{a=(u,v)\in\B}{u\mbox{ \tiny and }v\mbox{ \tiny differ in
      }j}}\j'(a,z)
\leq d\,N^{d-j+1}\,N^{j}=d\,N^{d+1},
\end{equation}
which eventually yields
\begin{equation}
\label{sc20}
\tau^\o_{21}\leq\mbox{const }\leps^{2d}\left(\max_{x\in S}\frac1{\o(x)}\right).
\end{equation}

\noindent{\bf Estimation of $\tau^\o_{22}$.} 
\begin{equation}
\label{sc21}
\tau^\o_{22}\leq\left(\max_{x\in S}\frac1{\o(x)}\right)
                N^{-d}\max_{b\in\B_2}\sum_{b'\in\B_2}\n(b,b').
\end{equation}
By Remark~\ref{good}, if $b,b'\in\B_2$ is in $\pi_{x,y}\in\Gamma$, then
$b'$ must be either in $Q_x$ or in $Q_y$ (see definition right
above Remark~\ref{good}). Thus for $b\in\B_2$, we have
\begin{equation}
\nn
\sum_{b'\in\B_2}\n(b,b')\leq\sum_{x,y\in S}\sum_{b'\in\B_2}
 [1\{b,b'\in Q_x\}+1\{b,b'\in Q_y\}
+1\{b\in Q_x,b'\in Q_y\}+1\{b\in Q_y,b'\in Q_x\}].
\end{equation}
Now
\begin{equation}
\label{sc23}
\sum_{x,y\in S}\sum_{b'\in\B_2}1\{b,b'\in Q_x\}
\leq
\sum_{y\in S}\sum_{x\in S}1\{b\in Q_x\}\sum_{b'\in\B_2}1\{b'\in Q_x\}.
\end{equation}
The two inner summands in the left hand side of~(\ref{sc23}) are uniformly 
bounded by const $\leps^{d}$, so the left hand side of~(\ref{sc23}) is bounded
by const $\leps^{2d}\,N^d$. For similar reasons, the same bound holds for 
$\sum_{x,y\in S}\sum_{b'\in\B_2}1\{b,b'\in Q_y\}$, 
$\sum_{x,y\in S}\sum_{b'\in\B_2}1\{b\in Q_x,b'\in Q_y\}$ and
$\sum_{x,y\in S}\sum_{b'\in\B_2}1\{b\in Q_y,b'\in Q_x\}$, and thus,
from~(\ref{sc21})
\begin{equation}
\label{sc24}
\tau^\o_{22}\leq\mbox{const }\leps^{2d}
\left(\max_{x\in S}\frac1{\o(x)}\right).
\end{equation}

The result of Proposition~\ref{prop:bound} now follows 
from~(\ref{sc6}),~(\ref{sc13}),~(\ref{sc20}),~(\ref{sc24}) and~(\ref{sc1}).
$\Box$


\vskip 1cm
\section{Decay of the annealed return probability for random walks on $\Z^d$}
\setcounter{equation}{0}
\label{sec:dec}


We go back to the study  of Markov chains taking their values in $\Z^d$.  
Let $\o:\Z^d\to\R^*_+$, 
and define   the Markov generator
\begin{equation}
\label{gen'}
\L^\o f(x)=\sum_{y\sim x}[\o(x)\wedge\o(y)]\,[f(y)-f(x)],
\end{equation} 
where the sum is over sites $y$ which are nearest neighbors to $x$.

As in Section~\ref{sec:comp},  
$\{X_t,\,t\in\R_+\}$ will be the coordinate process on path space 
$(\Z^d)^{\R_+}$ and we use the notation  $\P_x^\o$ to denote the unique probability measure 
on path space under which $\{X_t,\,t\in\R_+\}$ is the Markov process 
generated by~(\ref{gen'}) and satisfying $X_0=x$. 


As in Section~\ref{sec:times}, 
we choose the family $\{\o(x),\,x\in\Z^d\}$ at random, 
according to a law $\Q$ on
$(\R^*_+)^{\Z^d}$ such that 
\beqn
\nn
\mbox{ the random variables }\!\!\!\!\! &\{&\!\!\!\!\!\o(x),\,x\in\Z^d\} \mbox{ are i.i.d. };\\ 
\nn
\o(x)&\!\!\!\leq\!\!\!&1 \mbox{ for all }x;\\ 
\label{eq:cond3}
\Q(\o(0)\leq a)&\!\!\!\sim\!\!\!& a^\g \mbox{ as }a\downarrow0,
\eeqn
where $\g>0$ is a parameter.

\brm
We note that this generator has the same form as 
${\cal G}^\o$ in~(\ref{gen0}) by making $\o(x,y)=\o(x)\wedge\o(y)$,
and also the same form as 
${\cal L}^{\o,N}$ in~(\ref{gen}), but in infinite volume.
There would also be similar results for $\o$ defined on 
edges, instead of points, 
with i.i.d.~values for different edges, and the same technique 
would apply.
\erm

\brm

If $\o(0)$ were a Bernoulli random variable, then we would have a random walk
on a (independent, site) percolation cluster (provided we started in an
infinite occupied cluster). See~\cite{kn:MR1}.

\erm

In the sequel $\Q.\P_x^\o$ will be used as a short hand notation for the annealed 
law defined by $\Q.\P_x^\o[\,\cdot\,]=\int P_x^\o[\,\cdot\,]d\Q(\o)$. 
We are interested in estimating the decay of the return probability under $\Q.\P^\o$, 
$\Q.\P^\o[X_t=0]$, as $t$ tends to $+\infty$.  
It is actually quite easy to derive lower bounds for $\Q.\P^\o[X_t=0]$. 
Indeed, on one hand, one can use the comparison lemma~\ref{compalemma} with the usual nearest neighbor random 
walk on $\Z^d$ to prove that 
\beqn
\label{eq:low1}
\Q.\P^\o[X_t=0]\geq c t^{-d/2}, 
\eeqn 
for some contant $c$ that depends on the dimension $d$. 
There is another way to prove (\ref{eq:low1}), as follows. 
It is known~\cite{kn:DFGW} 
that, under $\Q.\P^\o$, $X_t$ satisfies the central 
limit theorem. Together with the reversibility and the translation invariance of 
the law $\Q$, the C.L.T. implies (\ref{eq:low1}) 
(See Appendix D, in~\cite{kn:MR1}). 

On the other hand, 
for any realization of $\o$, the first jump of $X_t$ follows 
an exponential law of parameter $\sum_{y\sim 0}\o(0)\wedge\o(y)\leq 2 d \o(0)$. 
Therefore
\beq \nn
\P^\o_0[X_t=0] \geq \P^\o_0[X_s=0, \forall s\leq t]= 
e^{ -t \sum_{y\sim 0}\o(0)\wedge\o(y)} 
\geq e^{-2 d \o(0) t}.
\eeq
Taking expectation w.r.t. $\Q$ and using the condition (\ref{eq:cond3}) on the law of $\o(0)$, 
a simple computation leads to a lower bound of the form 
\beqn
\label{eq:low2}
\Q.\P^\o[X_t=0]\geq c t^{-\g}.
\eeqn

As is indicated in the next statement, these lower bounds turn out to be of 
the correct logarithmic order. 

\begin{theo}
\label{t}
\beqn
\label{up}
\lim_{t\rightarrow+\infty}\frac{\log\Q.\P^\o[X_t=0]}{\log t}=-\left(\frac d2\wedge\g\right). 
\eeqn
\end{theo}

\brm
{}From the point of view of statistical mechanics --- here the statistical mechanics 
of a disordered system --- we consider Theorem~\ref{t} as a (nice) example of a 
dynamical phase transition. 
\erm

\brm
Such tools as Sobolev embeddings, isoperimetric or Nash inequalities of 
constant use for estimating transition probabilities of Markov chains, see~\cite{kn:Co},  
cannot be directly applied here because of the lack of ellipticity of the 
transition rates $\o$. Thus~(\ref{up}) is also an example of exotic 'heat kernel decay' 
for a non uniformly elliptic generator. 
\erm

\brm
A fruitful technique to handle r.w.r.e. is to isolate the effect of the fluctuations of the 
environment $\o$ in a given scale. See for instance random walks in Poisson environments~\cite{kn:Sz} where 
one single eigenvalue dominates the rest of the spectrum. There does not seem to 
exist such a separating scale in our model. 
\erm

In view of (\ref{eq:low1}) and (\ref{eq:low2}), only the upper bound is 
missing in the proof of (\ref{up}).

We use spectral theory. We rely on a trace formula similar to the one obtained
in Section~\ref{sec:comp}  
and on our spectral gap estimates from Proposition~\ref{prop:bound}. 


\subsection{Trace formula} 
\label{traceform}

We  express the annealed return probability as a trace. The argument is the same 
as in Section~\ref{sec:comp}, except that we restrict ourselves to computing the trace on cubes whose radius 
can be chosen as a function of time. This is possible because rates are assumed to be uniformly bounded.

Let $\xi>0$. In the sequel, we shall use the notation $N=t^{(1+\xi)/2}$. 
(In fact, $N$ should be defined as the integer part of $t^{(1+\xi)/2}$, but, 
for notational ease, we will omit integer parts.) 

Let $B_N=[-N,N]^d$, be the box centered at the origin and of radius $N$. Let 
$\L^{\o,N}$ be the restriction of the operator $\L^\o$ to $B_N$. 
Thus $\L^{\o,N}$ is defined by 
\begin{equation}
\L^{\o,N} f(x)=\sum_{y\sim x}[\o(x)\wedge\o(y)]\,[f(y)-f(x)],
\end{equation} 
where the sum is now restricted to neighboring points $x$ and $y$ in $B_N$ and we impose 
periodic boundary conditions. 
$-\L^{\o,N}$ is then a symmetric operator. 
We denote by $\{\l_i^\o(B_N), i\in[1,\# B_N]\}$ 
the set of its eigenvalues in increasing order. 

Let $\tau_N$ be the exit time of $X_t$ outside $B_N$. 

We compute $\Q.\P^\o_0[X_t=0]$ using the translation invariance of the probability $\Q$. 
Since $\Q.\P^\o_x[X_t=x]$ does not depend on $x$, we have 
\begin{eqnarray*}
\Q.\P^\o_0[X_t=0] &=& 
\frac 1{\#B_N} \sum_{x\in B_N} \Q.\P^\o_x[X_t=x]\\
&=& \frac 1{\#B_N} \sum_{x\in B_N} \Q.\P^\o_x[X_t=x;t<\tau_{2N}] 
    + \frac 1{\#B_N} \sum_{x\in B_N} \Q.\P^\o_x[X_t=x;t\geq \tau_{2N}]\\ 
&\leq&  \frac 1{\#B_N} \sum_{x\in B_{2N}} \Q.\P^\o_x[X_t=x;t<\tau_{2N}] 
       + \frac 1{\#B_N} \sum_{x\in B_N} \Q.\P^\o_x[t\geq \tau_{2N}].  
\end{eqnarray*}
If under $\P^\o_x$, $x\in B_N$, we have $t\geq \tau_{2N}$, then the process 
must have left the ball $x+B_N$ before time $t$. Since the probability 
$\Q.\P^\o_x[\exists s\leq t\mbox{ s.t. } X_s\notin x+B_N]$ does not depend on $x$, we 
have that $ \Q.\P^\o_x[t\geq \tau_{2N}]\leq \Q.\P^\o_0[t\geq \tau_N]$. 

We note that $ \sum_{x\in B_{2N}} \P^\o_x[X_t=x;t<\tau_{2N}] $ is the trace of 
the semi-group of the process $X_t$ killed when leaving the box $B_{2N}$, 
i.e., with Dirichlet boundary conditions outside $B_{2N}$. It is 
therefore dominated by the trace of $\exp(t\L^{\o,2N})$, that is  
\beqnn
\sum_{x\in B_{2N}} \P^\o_x[X_t=x;t<\tau_{2N}] 
\leq \sum_{i} e^{-\l_i^\o(B_{2N})t}.
\eeqnn

Thus, we have proved that 
\beqnn
\Q.\P^\o_0[X_t=0]                              
\leq \frac 1{\#B_N} \sum_{i} \Q[e^{-\l_i^\o(B_{2N})t}]+\Q.\P^\o[t\geq\tau_{N}].
\eeqnn

{}From the Carne-Varopoulos inequality, it follows that 
\beqn
\label{carvar}
\P^\o[t\geq\tau_{N}]
\leq 2tN^{d-1}e^{-\frac {N^2}{4t}}+e^{-ct}, 
\eeqn
where $c$ is a numerical constant, see Appendix C in~\cite{kn:MR1}.  
With our choice of $N=t^{(1+\xi)/2}$, we get that 
$\P^\o[t\geq\tau_N]$ decays faster than any polynomial as $t$ tends to $+\infty$. 

Thus Theorem \ref{t} will be proved if we can check that 
\beqn
\label{trace}
\lim_{\xi\rightarrow 0}\limsup_{t\rightarrow+\infty} 
\frac{\log \Q[\sum_{i} e^{-\l_i^\o(B_{N})t}]}{\log t}\leq 0\vee\left(\frac d 2-\gamma\right). 
\eeqn


\subsection{Min-Max} 

$C$ is a constant that depends only on $d$ and $\Q$. For constants depending on other 
parameters, we indicate it. 

Let us first recall the lower bound on the first non trivial eigenvalue of an 
operator of the form $\L^{\o,N}$. In Section~\ref{sec:times}, we proved that 
\beqn
\label{gap}
\frac 1 {\l^\o_2(B_N)}
\leq C\left(N^{2+\eps}+\sup_{x\in B_N} \frac 1{\o(x)}\right)d_\eps^N.
\eeqn  
In (\ref{gap}), $\eps$ is any positive number; 
$C$ is a constant depending on the dimension only; 
$d_\eps^N$ is a measure of the set 
$\{x\in B_N:\, \o(x)\leq N^{-\eps}\}$. 

With the notation of Section~\ref{sec:times}, Proposition~\ref{prop:bound}, 
$d_\eps^N=(\ell_\eps+1)^{2d}$. (But note that $\ell_\eps$ depends on $N$.) 
Thus $d_\eps^N$ is 
a random variable, i.e., depends on $\eps$, $N$ and also $\o$. 

Using the properties of $\Q$, we get that, for some constant $c$, that depends 
on $\Q$ only, we have 
\beqn
\label{despN}
\Q(d_\eps^N\geq A)
\leq c N^{-\frac{\eps\g A}2},  
\eeqn
where $A$ can be chosen  such that $A\geq\frac{4d}{\eps\g}$ and $N$ is supposed to be large enough. 
(How large depends on the dimension only.) A proof of (\ref{despN}) can be found 
in the proof of Lemma~\ref{lm:prob}. 

{}From the min-max caracterization of the eigenvalues of symmetric operators, 
we have 
\beqn
\nn
\l_{i+1}^\o(B_N)
=\max_{f_1,...,f_i}\min_{f}
\frac 12 \frac{ \sum_{x\sim y\in B_N} [\o(x)\wedge\o(y)]\, [f(x)-f(y)]^2}{\sum_{x\in B_N} f^2(x)},
\eeqn
where the 'max' is computed on choices of $i$ functions defined on $B_N$ and the 'min' 
is computed on functions $f$ such that, for all $j\in[1,i]$, $\sum_{x\in B_N} f(x)f_j(x)=0$. 

Thus, in the computation of $\l_{i+1}^\o(B_N)$, 
we may impose at most $i$ different linear constraints on the test function $f$. 
We consider two kind of conditions. 

Let $k\in\N^*$. We chop $\Z^d$ into a disjoint union of boxes of radius $k$, 
say $\Z^d=\cup_{z\in\Z^d} {\mathbf B}_z$, where 
${\mathbf B}_z=(2k+1)z+B_k$. 
We now choose for some of the function $f_j$'s, the indicator 
function of the boxes ${\mathbf B}_z$ that intersect $B_N$, i.e., we require that 
\beqnn
\sum_{x\in B_N\cap {\mathbf B}_z}f(x)=0, 
\eeqnn
for all $z\in\Z^d$ such that $B_N\cap{\mathbf B}_z\not=\emptyset$. The number of such $z$'s 
is at most 
\beqnn
n_2=\left(\frac{2N+1+2k+1}{2k+1}\right)^d.
\eeqnn 
Clearly, 
\beqnn
\sum_{x\in B_N}f^2(x)=\sum_z\sum_{x\in B_N\cap{\mathbf B}_z} f^2(x), 
\eeqnn
and
\beqnn
\sum_{x\sim y\in B_N} \o(x)\wedge\o(y) (f(x)-f(y))^2
\geq 
\sum_z\sum_{x\sim y\in B_N\cap{\mathbf B}_z} [\o(x)\wedge\o(y)]\,[f(x)-f(y)]^2.
\eeqnn
Therefore
\beqn 
\nn
\frac{ \sum_{x\sim y\in B_N} [\o(x)\wedge\o(y)]\,[f(x)-f(y)]^2}{\sum_{x\in B_N} f^2(x)} 
\geq
\min_z 
\frac{ \sum_{x\sim y\in B_N\cap{\mathbf B}_z} [\o(x)\wedge\o(y)]\,[f(x)-f(y)]^2}
{\sum_{x\in B_N\cap{\mathbf B}_z}  f^2(x)}, 
\eeqn
where, for each $z\in\Z^d$,   
$\sum_{x\in B_N\cap{\mathbf B}_z}f(x)=0$.

Next, let us choose $n_1$ points in $B_N$, say $\d_1,...,\d_{n_1}$. 
We choose for some of the $f_j$'s, the indicator function of the points 
$\d_j$ and their neighbors in $B_N$, i.e., we specify that 
$f(x)=0$, for $x\in\{\d_1,...,\d_{n_1}\}$ or $x\sim\d_j$, for some $j$. 
This recipe leads to, at most, $(2d+1)n_1$ different conditions. 
We note that, for such a function $f$, the value of the Dirichlet form 
$$\sum_{x\sim y\in B_N} [\o(x)\wedge\o(y)]\,[f(x)-f(y)]^2$$
does not depend on the value of $\o(\d_j)$ anymore. 
Therefore, we may assume that $\o(\d_j)=1$, for $j\in[1,n_1]$. 

Thus we see that, if $i\geq n_2+(2d+1)n_1$, then 
\beqnn
\l_{i+1}^\o(B_N)\geq\min_z \l_2^\oo(B_N\cap{\mathbf B}_z), 
\eeqnn
where $\oo$ is a new environment obtained by modifying the 
value of $\o$ to $1$ on all points $\d_i$ and $z$ ranges through those points in $\Z^d$ 
such that ${\mathbf B}_z$ intersects $B_N$. 

We now choose for $\d_j$ the points in $B_N$ where $\o$ achieves its lowest values. 
Let us use~(\ref{gap}) to estimate each eigenvalue 
$\l_2^\oo(B_N\cap{\mathbf B}_z)$: 
\beqn 
\label{eigenvalueestimate}
\frac 1 {\l_{i+1}^\o(B_N)}
\leq C\left(k^{2+\eps}+\sup_{x\in B_N}^{n_1}\frac 1{\o(x)}\right)d_\eps^N. 
\eeqn

We used $d_\eps^N$ as a uniform upper bound for the minimal side length of 
strips for which the event $A_N(L)$ in Definition~\ref{ell}  occurs. 

$\sup_{x\in B_N}^{n_1} 1/\o(x)$ denotes the maximal value 
of $1/\oo(x)$, i.e., 
\beqnn
\sup_{x\in B_N}^{n_1} \frac1{\o(x)}
=\max \{ h:\, \#\{x\in B_N:\, \o(x)=1/h\}\geq n_1+1\}.
\eeqnn


\subsection{Proof of Theorem \ref{t}} 

Remember that we have already chosen some parameter $\xi>0$ (that we want to 
choose close to $0$ and which is related to $t$ 
by $N=t^{(1+\xi)/2}$), and another parameter $\eps>0$ which is arbitrarily close to $0$. 
We need a third parameter  $a\in(0,1)$. 
The constant $A$ in~(\ref{despN}) is at our disposal. 
We also still have to choose $n_1$ and $n_2$, depending on $i$ and such that 
$i\geq n_2+(2d+1)n_1$. 

Write  
\beqn \nn
&&\Q\left[\sum_i e^{-\l_i^\o(B_N) t}\right]\\\nn
&\leq& \Q\left[\sum_i e^{-\l_i^\o(B_N) t};d^N_\eps\geq A\right]
+\sum_i\Q\left[e^{-\l_i^\o(B_N)t};\l_i^\o(B_N)\geq 
N^{-\eps}\left(\frac{i^{a/d}}N\right)^2\right]\\\nn
&+&\sum_i \Q\left[e^{-\l_i^\o(B_N)t}; d^N_\eps\leq A\mbox{ and } 
\l_i^\o(B_N)\leq N^{-\eps} \left(\frac {i^{a/d}}N\right)^2\right] \\\nn
&\leq &(2N+1)^d \Q[d^N_\eps\geq A] 
+\sum_i e^{-N^{-\eps} \left(\frac {i^{a/d}}N\right)^2 t}\\
\label{equation}
&+&\sum_i \Q\left[d^N_\eps\leq A\mbox{ and } 
\l_i^\o(B_N)\leq N^{-\eps} \left(\frac {i^{a/d}}N\right)^2\right]. 
\eeqn 
Using (\ref{despN}), we see that we can choose $A$ in such a way that 
\beqn 
\label{bound1}
\limsup_{t\rightarrow+\infty} \frac {\log [(2N+1)^d \Q(d^N_\eps\geq A)]}
{\log t}\leq 0\vee\left(\frac d 2-\g\right).
\eeqn

An easy computation shows that 
\beqn 
\label{bound2}
\limsup_{t\rightarrow+\infty} 
\frac {\log \sum_i e^{-N^{-\eps}\left(\frac {i^{a/d}}N\right)^2 t}}{\log t}
\leq \frac d{2a}\left(\xi+\frac \eps 2+\frac{\eps\xi}2\right).
\eeqn

Let us now bound the last term in (\ref{equation}). 
Assume that $d^N_\eps\leq A$ and 
$\l_i^\o(B_N)\leq N^{-\eps} \left(\frac {i^{a/d}}N\right)^2$.
 From~(\ref{eigenvalueestimate}), we must have 
\beqnn
N^{\eps+2}i^{-\frac {2a}d}\leq C\left(k^{2+\eps}+\sup^{n_1}_{x\in B_N} \frac 1{\o(x)}\right)A.
\eeqnn

We choose $n_2=i^a$ and assume that $i$ is large enough, how large depending on the 
dimension, $a$ and $\g$ only, which we may do. 
Then $k^{2+\eps}\leq N^{2+\eps} i^{-\frac{a(2+\eps)}d}$. Therefore, we must have 
\beqnn
N^{\eps+2}i^{-\frac {2a}d}\leq C\sup^{n_1}_{x\in B_N} \frac 1{\o(x)}, 
\eeqnn
with a possibly different value for $C$. 

{}From now on, we deal separately with the cases of large or small values of $\g$. 

{\bf Case $\g\geq\frac d 2$}. We then choose $\eps<2\frac d\g$ and 
$a=1-\frac d{2\g}+\frac \eps 4$.

The computation goes as follows (the value of $C$ changes from line to line) 
\beqn\nn
&&\sum_i\Q\left[d_\eps^N\le A \mbox{ and } 
\l_i^\o(B_N)\leq N^{-\eps}\left(\frac{i^{a/d}}N\right)^2\right]\leq 
\sum_i \Q\left[C\sup^{n_1}_{x\in B_N} \frac 1{\o(x)}\geq N^{\eps+2}
i^{-\frac{2a}d}\right]\\ \nn
&\leq& \sum_i \Q\left[C\sup^{n_1}_{x\in B_N} \frac 1{\o(x)}\geq
  N^{\eps+2-2a}\right]= \Q\left[\#\left\{i:\, C\sup^{n_1}_{x\in B_N} 
\frac 1{\o(x)}\geq N^{\eps+2-2a}\right\}\right]\\ \nn
&\leq& C\Q\left[\#\left\{i:\, C\sup^i_{x\in B_N} \frac 1{\o(x)}\geq
    N^{\eps+2-2a}\right\}\right]= C (2N+1)^d  
\Q\left[\frac 1{\o(x)}\geq N^{\eps+2-2a}\right]\\ 
\label{bound0}
&\leq& C N^d N^{-\g(\eps+2-2a)}= C N^{-\frac{3\eps\g}4}, 
\eeqn 
where the second inequality follows because $i\leq (2N+1)^d$, and the
third one because $n_1= (i-i^a)/(2d+1)$ and $a<1$.
Thus we deduce from (\ref{bound0}), (\ref{equation}), (\ref{bound1}) 
and (\ref{bound2}) that 
\beqnn
\limsup_{t\rightarrow+\infty}\frac{\log \Q\left[\sum_i e^{-\l_i^\o(B_N)t}\right]}{\log t}
\leq \frac d{2a}\left(\xi+\frac\eps 2+\frac{\eps\xi}2\right).
\eeqnn
Let $\eps$ tend to $0$ and then $\xi$ tend to $0$ to deduce (\ref{trace}). 
This ends the proof of Theorem \ref{t} in the case $\g\geq\frac d 2$.

{\bf Case $\g<\frac d 2$}. Let $\d\in(0,\g)$, to be chosen later.
We have  
\beqn \nn
&\sum_i\Q\left[d_\eps^N\le A \mbox{ and } \l_i^\o(B_N)\leq N^{-\eps} 
\left(\frac{i^{a/d}}N\right)^2\right]\leq\sum_i \Q\left[C\sup^{n_1}_{x\in B_N} \frac 1{\o(x)}\geq N^{\eps+2}i^{-\frac{2a}d}\right]&\\ 
\nn
&\leq  N^{(2a-\eps-2)\d} \sum_i \Q\left[\left(\sup^{n_1}_{x\in B_N} \frac 1{\o(x)} \right)^\d\right],& 
\eeqn 
since $i\leq (2N+1)^d$. 
Remember that $n_2=i^a$ is much smaller than $n_1= (i-i^a)/(2d+1)$ for large values of $i$, say 
$i/(4d+2)\leq n_1 \leq i/(2d+1)$. Let 
$x_0,...x_j,...x_{(2N+1)^d-1}$ be an enumeration of the points in $B_N$ such that the sequence 
$\o(x_j)$ is increasing. Thus 
\beqnn 
\sup^{n_1}_{x\in B_N} \frac 1{\o(x)}=\frac 1{\o(x_{n_1})}
\leq \frac 1{\o\left(x_{i/(4d+2)}\right)}. 
\eeqnn
Therefore 
\beqnn
\sum_i \Q\left[\left(\sup^{n_1}_{x\in B_N} \frac 1{\o(x)} \right)^\d\right] 
\leq (4d+2) \sum_{x\in B_N} \Q\left[\left(\frac 1{\o(x)}\right)^\d\right]
= c_\delta (4d+2)(2N+1)^d,
\eeqnn
where $c_\d=\Q\left[\left(\frac 1{\o(x)}\right)^\d\right]$. Note that $\Q\left[\left(\frac 1{\o(x)}\right)^\d\right]$ is finite and does not depend 
on $x$. Therefore 
\beqnn 
\sum_i\Q\left[d_\eps^N\le A \mbox{ and } \l_i^\o(B_N)\leq N^{-\eps} \left(\frac{i^{a/d}}N\right)^2\right]
\leq c_\d N^{(2a-\eps-2)\d} (2N+1)^d.
\eeqnn
Gathering this last inequality with (\ref{bound1}) and (\ref{bound2}), we get that 
\beqn 
\nn
&&\limsup_{t\rightarrow+\infty}\frac{\log \Q\left[\sum_i e^{-\l_i^\o(B_N)t}\right]}{\log t}\\ 
\label{c1}
&\leq &\max\left\{\frac d 2 -\g; \frac d{2a}\left(\xi+\frac\eps 2+\frac{\eps\xi}2\right) ; 
\frac{1+\xi}2[d+(2a-\eps-2)\d]\right\}, 
\eeqn
with $a\in(0,1)$ and $\d\in(0,\g)$. First replace $\d$ by $\g$. 
Then let $\eps$ tend to $0$ and choose $a=\frac {d\xi}{d-2\g}$. 
The upper bound in (\ref{c1}) becomes
$\max\left[\frac d 2 -\g; \frac{1+\xi}2(d-2\g)+\frac{(1+\xi)\xi d\g}{d-2\g}\right]$. 
Finally let $\xi$ tend to $0$ and conclude that 
\beqnn
\lim_{\xi\rightarrow 0}\limsup_{t\rightarrow+\infty}\frac{\log \Q\left[\sum_i e^{-\l_i^\o(B_N)t}\right]}{\log t} 
\leq \frac d 2 -\g,  
\eeqnn
and Theorem \ref{t} is now proved in the case $\g<\frac d 2$. $\Box$



\vskip 1cm
\section{Quenched decay of the return probability}
\setcounter{equation}{0}
\label{sec:quenc}


In this section, we investigate the quenched decay of the return probability. 
Model and notation are the same as in Section~\ref{sec:dec}: a random walk 
among i.i.d.~random conductancies with a power law with an exponent $\gamma$. 
Now we are rather interested in the asymptotics of the return probability 
$\P^\o_0[X_t=0]$ in $\Q$ probability. Let us set $\a_c$ to be the best exponent $\a$ 
such that 
\beqn \Q[\P^\o_0[X_t=0]\leq t^{-\a}]\rightarrow 1 
\hbox{ as } t\rightarrow \infty\,. \eeqn 
 From Theorem \ref{t}, it is clear that 
$\a_c\geq \frac d2\wedge\g$. We can do better in the case $\g<\frac d2$: 

\begin{theo}
\label{t'}
For any $\g<\frac d2$ then $\a_c>\g$. 
\end{theo}

\brm 
Although rather unsatisfactory --- because it does not give the true vaalue of $\a_c$ --- 
Theorem~\ref{t'} 
shows that the typical decay of the return probability 
is strictly faster than the averaged decay. Such a situation is sometimes called in the litterature a 
'high disorder regime'. 
\erm

\brm 
The proof of Theorem~\ref{t'} actually yields the lower bound 
\beqn 
\a_c\geq \frac d 2 \frac {1+\g}{1+d/2} . \eeqn 
There is no reason to believe that this bound is sharp for a given value of $\g$. 
Notice however that, in the regime $\gamma\rightarrow \frac d 2$, we get the inequality 
$\a_c\geq \frac d 2$, which seems to be sharp. 
\erm 

Let us sketch the proof: we use the fact that, with large $\Q$ probability, 
the origin lies in an infinite percolation cluster, say $\C$, of 'good' sites, 
where $\o$ is bounded from below. 
Estimates on the return probability for random walks on percolation clusters
have been proved in~\cite{kn:MR1} (See also~\cite{kn:Ba} ). One strategy would then be to try to couple 
the random walk in the environment $\o$ with the random walk on $\C$: 
we have no idea on how to do that. 
We rather rely on spectral theory to compare the behaviours of the eigenvectors for the two random 
walks. 
Note that from the results of~\cite{kn:MR1} follow precise estimates on the eigenvalues of the 
discrete 
Laplace operator on $\C$. The core of the proof is to show that eigenvectors of the generator of the 
random walk in the environment $\o$, when they correspond to small enough eigenvalues, are 
concentrated 
outside $\C$, and therefore do not contribute too much to the asymptotics of the return probability as 
soon as the random walk starts outside $\C$. 


\subsection{Proof of Theorem~\ref{t'}. Step 1}

Let $\a<\frac d 2 \frac {1+\g}{1+d/2}$. 
Choose two parameters $\eps>0$ and $\xi>0$.   We shall use the notation $N=t^{(1+\xi)/2}$. 
(In fact, $N$ should be defined as the integer part of $t^{(1+\xi)/2}$, but, 
for notational ease, we will omit integer parts.)  All the limits to be taken 
are to be understood as $t\rightarrow\infty$ or, equivalently $N\rightarrow\infty$. 

Let $\C^\o$ be the largest connected component of the set 
$\{x\in\Z^d\, : \, \o(x)\geq N^{-\eps}\}$. 
We assume that $N$ is large enough so that 
$\Q[\o(x)\geq N^{-\eps}]$ becomes larger than the critical percolation probability 
on $\Z^d$. Then  $\C^\o$ is the unique infinite connected component of 
the set $\{x\in\Z^d\, : \, \o(x)\geq N^{-\eps}\}$, see~\cite{kn:G}. 
We denote by $\C^\o_N$ the largest connected component of the intersection 
$\C^\o\cap B_N$, where $B_N=[-N,N]^d$. 

In the next step of the proof, we will define a set of environments, denoted $\O_N$, such that 
$\Q[\O_N]\rightarrow 1$. We further have the property 
$\Q[\frac{\#\C^\o_N}{\#B_N}]\rightarrow 1$. 

Calling $\{\l_i^\o(B_N), i\in[1,\# B_N]\}$ the eigenvalues of $-\L^{\o,N}$ in increasing order, 
and $\{\psi^\o_i, i\in[1,\# B_N]\}$ the corresponding eigenvectors with due normalization in 
$L^2(B_N)$, a very similar computation as in Subsection~\ref{traceform} leads to the following 
series of inequalities.  

We first use the invariance by translation of $\Q$. 
\beqnn 
\Q[\P^\o_0[X_t=0]\geq t^{-\a}]
=\Q[\P^\o_x[X_t=x]\geq t^{-\a}] \eeqnn 
holds for any $x\in B_N$. Therefore 
\beqnn \Q[\P^\o_0[X_t=0]\geq t^{-\a}] 
=\frac 1{\#B_N} \sum_{x\in B_N} \Q[\P^\o_x[X_t=x]\geq t^{-\a}] .\eeqnn 

Note that 
\beqnn  
\Q[\P^\o_x[X_t=x]\geq t^{-\a}]
\leq \Q[\P^\o_x[X_t=x; t<\tau_{2N}]\geq t^{-\a}/2]  
+ \Q[\P^\o_x[t\geq\tau_{2N}]\geq t^{-\a}/2]\,, \eeqnn   
where $\tau_{2N}$ is the exit time of $B_{2N}$. Since 
$\sup_\o\sup_x\P^\o_x[t\geq\tau_{2N}]$ decays faster than any polynomial, see~(\ref{carvar}), we have 
\beqnn \limsup \Q[\P^\o_0[X_t=0]\geq t^{-\a}] 
\leq \limsup \frac 1{\#B_N} \sum_{x\in B_N} \Q[\P^\o_x[X_t=x; t<\tau_{2N}]\geq t^{-\a}/2].\eeqnn 

We now restrict our attention to those environments belonging to $\O_N$ and to the points 
$x\in\C^\o_N$: 
\beqnn  
&&\Q[\P^\o_x[X_t=x; t<\tau_{2N}]\geq t^{-\a}/2]\\  
&\leq& \Q[\P^\o_x[X_t=x; t<\tau_{2N}]\geq t^{-\a}/2; x\in\C^\o_N; \O_N]
+\Q[\O_N^c]+\Q[x\notin\C^\o_N].\eeqnn 

Since $\Q[\O_N^c]\rightarrow 0$, we therefore get that 
\beqnn 
&&\limsup \Q[\P^\o_0[X_t=0]\geq t^{-\a}]\\  
&\leq& \limsup \frac 1{\#B_N} \sum_{x\in B_N} 
  \Q[\P^\o_x[X_t=x; t<\tau_{2N}]\geq t^{-\a}/2; x\in\C^\o_N; \O_N] 
  +\limsup \Q\left[\frac{\#(\C^\o_N)^c}{\#B_N}\right].\eeqnn
But since $\Q[\frac{\#\C^\o_N}{\#B_N}]\rightarrow 1$ (see step 2 below), we have
\beqnn 
&&\limsup \Q[\P^\o_0[X_t=0]\geq t^{-\a}] \\ 
&\leq &\limsup \frac 1{\#B_N} \sum_{x\in B_N} 
  \Q[\P^\o_x[X_t=x; t<\tau_{2N}]\geq t^{-\a}/2; x\in\C^\o_N; \O_N].\eeqnn 

 From the Markov inequality, we deduce that 
\beqnn 
\Q[\P^\o_x[X_t=x; t<\tau_{2N}]\geq t^{-\a}/2; x\in\C^\o_N; \O_N]
\leq 2t^\a \Q[\P^\o_x[X_t=x; t<\tau_{2N}]; x\in\C^\o_N; \O_N],\eeqnn 
and thus 
\beqnn 
\limsup \Q[\P^\o_0[X_t=0]\geq t^{-\a}]  
\leq 2 \limsup \frac {t^\a} {\#B_N} 
\Q[ \sum_{x\in \C^\o_N} \P^\o_x[X_t=x; t<\tau_{2N}]; \O_N].\eeqnn 

Finally we express the probability $\P^\o_x[X_t=x; t<\tau_{2N}]$ in the spectral decomposition 
as 
\beqnn \P^\o_x[X_t=x; t<\tau_{2N}]
=\frac 1{\#B_N} \sum_i e^{-\l_i^\o(B_{2N})\, t} (\psi^\o_i(x))^2,\eeqnn 
and get that 
\beqn\label{eq:4} 
\limsup \Q[\P^\o_0[X_t=0]\geq t^{-\a}]  
\leq 2 \limsup \frac {t^\a} {\#B_N} 
\Q\left[\sum_i e^{-\l_i^\o(B_{2N})\,t} \frac 1{\#B_N} \sum_{x\in \C^\o_N}(\psi^\o_i(x))^2; \O_N\right].\eeqn 

Let us pause a little to look at (\ref{eq:4}). It is true that 
$\frac 1{\#B_N} \sum_{x\in \C^\o_N}(\psi^\o_i(x))^2\leq 1$; but if we would use this upper bound, 
we would be left with $\Q[\sum_i e^{-\l_i^\o(B_{2N})\,t}]$, and the best value for $\a$ would then be 
$\gamma$, as the results of Section~\ref{sec:dec} show. We have to find a better way. 
Note that terms corresponding to large values of $i$, and thus 
large values of $\l_i^\o(B_{2N})$, can be easily controlled. Thus the main point is to 
show that $\frac 1{\#B_N} \sum_{x\in \C^\o_N}(\psi^\o_i(x))^2$ is small enough for 
small $i$, i.e. we have to prove that eigenvectors corresponding to small eigenvalues are 
concentrated outside $\O_N$. And in fact one would expect this to be true since small eigenvalues 
arise because of small values of $\o$, and these precisely sit outside $\C^\o_N$.


\subsection{Step 2. Definition of $\O_N$} 

The set $\O_N$ is defined by two requirements: we ask that for any $\o\in\O_N$ we have 
$$\mbox{(i)}\,\, 0\in\C^\o_N.$$

The second requirement deals with the behaviour of the random walk on   $\C^\o_N$: 
let $(\mu_i, i\in[1,\#\C^\o_N])$ be the eigenvalues of the discrete Laplace operator on $\C^\o_N$ as defined 
in~\cite{kn:MR1}. We will also use the notation $(\phi_i, i\in[1,\#\C^\o_N])$ for the corresponding eigenvectors. 
We assume that the eigenvalues are in increasing order and the eigenvectors are normalized in 
$L^2(\C^\o_N)$ for the counting measure. Of course the $\mu_i$s and $\phi_i$s depend on $\o$ and 
$N$. 

Let \beqnn \eta=\frac{1+\g}{\frac 12+\frac 1d}\, \hbox{ and }\,  j=N^{d-\eta}.\eeqnn 
Note that since $\g<d/2$, then $\eta<d$. 
We then require that, on $\O_N$, 
$$\mbox{(ii)}\,\, \mu_j\geq \frac {j^{2/d}}{N^2(\log N)^{8(d-\eta)/d}}.$$

The definition of $\O_N$ is now complete and all that remains to be done is to
check that\\ $\Q(\O_N)\rightarrow 1$. 

That $\Q(\mbox{(i) holds})=\Q(0\in\C^\o_N)\rightarrow 1$ is obvious. 

As for condition (ii), we rely on the results of~\cite{kn:MR1}. Calling 
$P^\o_x[X^N_s=y]$ the transition probabilities for the random walk on $\C^\o_N$, we quote 
from formula (6) of \cite{kn:MR1}: $\Q$-a.s.~on the set where $\C^\o$ is infinite 
\beqnn 
\sup_{x,y\in\C^\o_N} \left\vert \frac 1{\#\C^\o_N}- P^\o_x[X^N_s=y]\right\vert 
\leq C \frac {(\log N)^{2d}}{s^{\frac d 2+d\frac{\log\log N}{\log N}}},\eeqnn 
where $C$ is a dimension dependent constant, $s$ is arbitrary, and $N\geq N_0(\o)$ is large enough. 
(In~\cite{kn:MR1}, formula (6) is deduced from the isoperimetric inequality (4), (4) is 
a consequence of (21), and (21) is proved for both site and bond percolation models with parameter 
$p$ close enough to $1$, which is our case here. Besides, we replaced $\eps(N)$ by its value  
$\eps(N)=d+2d\frac{\log\log N}{\log N}$, noticing that 
$(4\eps(N)/\beta^2)^{\eps(N)/2}$ then behaves like a constant.)  

We then choose $x=y$, sum over $x\in\C^\o_N$, and express the result as a trace to get that 
\beqnn 
\sum_i e^{-\mu_i s}\leq 1+\#\C^\o_N C  
\frac {(\log N)^{2d}}{s^{\frac d 2+d\frac{\log\log N}{\log N}}}.\eeqnn 
Therefore 
\beqnn 
j e^{-\mu_j s}\leq 1+C  \#\C^\o_N\frac {(\log N)^{2d}}{s^{\frac d 2+d\frac{\log\log N}{\log N}}}.\eeqnn 
Take now 
$s=\frac{N^2}{j^{2/d}}(\log N)^{8(d-\eta)/d}$. Then 
\beqnn 
\frac 1 j C \#\C^\o_N\frac {(\log N)^{2d}}{s^{\frac d 2+d\frac{\log\log N}{\log N}}}\rightarrow 0,
\eeqnn
so that $e^{-\mu_j s}\rightarrow 0$, and  
we have proved that $\Q$-a.s.~on the set where $\C^\o$ is infinite, for large enough $N$, 
condition (ii) is fullfilled. 

Finally we already used the fact that $\Q[\frac{\#\C^\o_N}{\#B_N}]\rightarrow 1$ that should be justified: 
from the result of Appendix B 
of~\cite{kn:MR1}, we know that the expected density in $B_N$ of the component of
$\C^\o\cap B_N$ that contains 
the origin goes to 1 as $N\to\infty$, and Lemma~\ref{etab} implies that the 
expected density in $B_N$ of the largest component of $\C^\o\cap B_N$  
goes to 1 as $N\to\infty$.
Thus the  component of
$\C^\o\cap B_N$ that contains 
the origin and $\C^\o_N$ coincide for large $N$ and 
its density tends to $1$.


\subsection{Step 3. Spectral analysis} 

Assume that $\o\in\O_N$. 

We bound the term $\sum_{x\in B_N}(\psi^\o_i(x))^2$ in (\ref{eq:4}) in two steps 
by writing that 
\beqnn \frac 1{\#B_N} \sum_{x\in \C^\o_N}(\psi^\o_i(x))^2 
= \frac 1{\#B_N} \sum_{x\in \C^\o_N}(\psi^\o_i(x)-P^j\psi^\o_i(x))^2+ 
\frac 1{\#B_N} \sum_{x\in \C^\o_N}(P^j\psi^\o_i(x))^2,\eeqnn 
where $P^j$ is the projection on the subspace of $L^2(\C^\o_N)$ spanned by the 
eigenvectors $(\phi_i,i\in[1,j])$. 

On one hand, since $ \frac 1{\#B_N} \sum_{x\in B_N} (\phi_i(x))^2\leq 1$, then 
\beqnn  
&&\sum_i e^{-\l_i^\o(B_{2N})\,t} \frac 1{\#B_N} \sum_{x\in \C^\o_N}(P^j\psi^\o_i(x))^2\\ 
&\leq& \sum_i \frac 1{\#B_N} \sum_{x\in \C^\o_N}(P^j\psi^\o_i(x))^2\\
&=& \sum_i\sum_{k\leq j} \frac 1{\#B_N} \sum_{x\in \C^\o_N}\phi_k(x)\psi^\o_i(x)\\ 
&=& \sum_{k\leq j} \frac 1{\#B_N} \sum_{x\in \C^\o_N}(\phi_k(x))^2\\ 
&\leq& j=N^{d-\eta}.\eeqnn 

On the other hand, for any function $f$ on $\C^\o_N$, we have 
\beqnn \frac 1{\#\C^\o_N} \sum_{x\in \C^\o_N}(f(x)-P^jf(x))^2
\leq \frac 1{\mu_j} \frac 1{2\#\C^\o_N}  \sum_{x\sim y\in \C^\o_N}(f(x)-f(y))^2 , 
\eeqnn this last expression being the Dirichlet form of the random walk on $\C^\o_N$. 
Since $\o(x)\geq N^{-\eps}$ on $\C^\o_N$, we get 
\beqnn \sum_{x\sim y\in \C^\o_N}(f(x)-f(y))^2 
\leq N^\eps \sum_{x\sim y\in B_{2N}}(\o(x)\wedge\o(y))(f(x)-f(y))^2,\eeqnn 
this last expression being now the Dirichlet form of the random walk on $B_{2N}$. 
Since $\psi^\o_i$ is an eigenvector, 
\beqnn 
\frac 1{2\#B_{2N}} \sum_{x\sim y\in B_{2N}}(\o(x)\wedge\o(y))(\psi^\o_i(x)-\psi^\o_i(y))^2
=\l_i^\o(B_{2N}).\eeqnn 
So  
\beqnn 
\frac 1{\#B_N} \sum_{x\in \C^\o_N}(\psi^\o_i(x)-P^j\psi^\o_i(x))^2 
\leq 2^dN^\eps \frac {\l_i^\o(B_{2N})}{\mu_j}.\eeqnn 

 From these two estimates, we deduce that 
\beqn\label{eq:4.2} 
\sum_i e^{-\l_i^\o(B_{2N})\,t} \frac 1{\#B_N} \sum_{x\in \C^\o_N}(\psi^\o_i(x))^2  
\leq N^{d-\eta} +2^d \frac{N^\eps}{\mu_j} \sum_i \l_i^\o(B_{2N})e^{-\l_i^\o(B_{2N})\,t}. 
\eeqn


\subsection{Step 4} 

Combining (\ref{eq:4}) and (\ref{eq:4.2}), we see that Theorem (\ref{t'}) will be proved once 
we have checked that 
$\frac {t^\a} {\#B_N} N^{d-\eta}\rightarrow 0$ and that 
\beqn\label{eq:4.3}  
\frac {t^\a} {\#B_N} {N^\eps} 
\Q\left[\frac 1{\mu_j}\sum_i \l_i^\o(B_{2N})e^{-\l_i^\o(B_{2N})\,t};\O_N\right]\rightarrow 0.\eeqn 

We recall that $\a<\frac d 2 \frac {1+\g}{1+d/2}$, $N=t^{(1+\xi)/2}$, 
$\eta=\frac{1+\g}{\frac 12+\frac 1d}$, and 
$\mu_j\geq (\frac {j^{1/d}}N)^2(\log N)^{-8(d-\eta)/d}$ on $\O_N$. 
It is then immediate to see that $\frac {t^\a} {\#B_N} N^{d-\eta}\rightarrow 0$. 
Besides, (\ref{eq:4.3}) will hold for any $\a<\frac d 2 \frac {1+\g}{1+d/2}$ 
and some $\eps>0$ if  
\beqn\label{5}
\lim_{\xi\rightarrow 0} \limsup \frac{ \Q[\sum_i \l_i^\o(B_{2N})e^{-\l_i^\o(B_{2N})\,t}]}{\log t} 
\leq\frac d 2 -1-\g.\eeqn

But using the inequality $ \l_ie^{-\l_it}\leq \frac 1 t e^{-\frac 1 2 \l_it}$, we get 
\beqnn 
&&\lim_{\xi\rightarrow 0} \limsup \frac{ \Q[\sum_i \l_i^\o(B_{2N})e^{-\l_i^\o(B_{2N})\,t}]}{\log t} \\ 
&\leq& -1+\lim_{\xi\rightarrow 0} \limsup \frac{ \Q[\sum_i e^{-\l_i^\o(B_{2N})\,t/2}]}{\log t}\\ 
&\leq& -1+\frac d 2 -\g,\eeqnn 
by~(\ref{up}). $\Box$

\end{document}